%% file: AllenCahn_131021.tex
%
\newif\ifloadreferences\loadreferencestrue
%
\input preamble %
%
%
%
%
%
\makeop{in}%
\makeop{ds}%
\makeop{Ev}%
\makeop{dt}%
\makeop{Lim}%
\makeop{Index}%
\makeop{AC}%
\makeop{pAC}%
\makeop{ell}%
\makeop{par}%
\makeop{sgn}%
\makeop{Sing}%
\makeop{Crit}%
\makeop{Eval}%
\makeop{Diag}%
\makeop{Null}%
\makeop{HAC}%
\newref{Nirenberg}{Agmon S., Nirenberg L., Lower bounds and uniqueness theorems for solutions of differential equations in a Hilbert space, {\sl Comm. Pure Appl. Math.}, {\bf 20}, (1967), 207--229}%
\newref{AllenCahn}{Allen S., Cahn J., A microscopic theory for antiphase boundary motion and its application to antiphase domain coarsening, {\it Acta Metall.}, {\bf 27}, (1979), 1084--1095}
\newref{Aronszajn}{Aronszajn N., A unique continuation theorem for solutions of elliptic partial differential equations or inequalities of second order, {\sl J. Math. Pures Appl.}, {\bf 36}, (1957), 235--249}%
\newref{Cabre}{Cabr\'e X., Uniqueness and stability of saddle-shaped solutions to the Allen-Cahn equation, {\sl J. Math. Pures Appl.}, {\bf 98}, (2012), no. 3, 239--256}
\newref{DeGiorgi}{De Giorgi E., Some conjectures on flow by mean curvature, {\sl White Paper}, (1990)}%
\newref{GilbTrud}{Gilbarg D., Trudinger N. S., {\sl Elliptic partial differential equations of second order}, Classics in Mathematics, Springer-Verlag, Berlin, (2001)}%
\newref{Ilmanen}{Ilmanen T., Convergence of the Allen-Cahn Equation to Brakke's motion by Mean Curvature, {\sl J. Diff. Geom.}, {\bf 38}, (1993), 417--461}
\newref{Kato}{Kato T., {\sl Perturbation theory for linear operators}, Classics in Mathematics, Springer-Verlag, Berlin, (1995)}%
\newref{Palais}{Palais R. S., Morse theory on Hilbert manifolds, {\sl Topology}, {\bf 2}, (1963), 299--340}
\newref{Salamon}{Robbin J., Salamon D., The spectral flow and the Maslov index, {\sl Bull. London Math. Soc.}, {\bf 27} (1995), no. 1, 1--33}%
\newref{SalamonWeber}{Salamon D., Weber J., Floer homology \& the heat flow, {\sl GAFA}, {\bf 16}, (2006), no. 5, 1050--1138}
\newref{Schwarz}{Schwarz M., {\sl Morse homology}, Progress in Mathematics, {\bf 111}, Birkh\"auser Verlag, Basel, (1993)}%
\newref{Smale}{Smale S., An infinite dimensional version of Sard's theorem, {\sl Amer. J. Math.}, {\bf 87}, (1965), 861--866}%
\newref{SmiMH}{Smith G., A H\"older Space Approach to Morse Homology, in preparation}%
\newref{Taylor}{Taylor M. E., {\sl Partial differential equations III. Nonlinear equations.}, Applied Mathematical Sciences, {\bf 117}, Springer, New York, (2011)}%
\newref{Weber}{Weber J., Morse homology for the heat flow, {\sl Math. Z.}, {\bf 275}, (2013), no. 1, 1--54}%
\myfontdefault
\global\headno=0
\global\showpagenumflag=1
\def\Pagetitle{}
\def\Pagefooter{}
\null
\vfill
\def\centre{\rightskip=0pt plus 1fil \leftskip=0pt plus 1fil \spaceskip=.3333em \xspaceskip=.5em \parfillskip=0em \parindent=0em}%
\def\textmonth#1{\ifcase#1\or January\or Febuary\or March\or April\or May\or June\or July\or August\or September\or October\or November\or December\fi}
\font\abstracttitlefont=cmr10 at 14pt
{\abstracttitlefont\centre Bifurcation of Solutions to the Allen-Cahn Equation\par}
\bigskip
{\centre Graham Smith\par}
\bigskip
{\centre \the\day\ \textmonth\month\ \the\year\par}
\bigskip
{\centre Instituto de Matem\'atica, UFRJ,\par
Av. Athos da Silveira Ramos 149,\par
Centro de Tecnologia - Bloco C,\par
Cidade Universit\'aria - Ilha do Fund\~ao,\par
Caixa Postal 68530, 21941-909 Rio de Janeiro,\par
RJ - BRASIL\par}
\bigskip
\noindent{\bf Abstract:\ }We use Morse Homology to study bifurcation of the solution sets of the Allen-Cahn Equation.
%
%
\par
\nextoddpage
\global\showpagenumflag=0
\global\pageno=1
\def\Pagetitle{\hfil The Allen-Cahn Equation\hfil}
\def\Pagefooter{\hfil{\myfontdefault\folio}\hfil}
\newhead{Introduction}
\newsubhead{Main Results} Let $M:=(M^n,g)$ be a compact, $n$-dimensional Riemannian manifold. For $\epsilon>0$, we define the {\bf Allen-Cahn Operator} over $C^\infty(M)$ with parameter $\epsilon$ by:
$$
\opAC_{\epsilon,g}(u) := -\epsilon\Delta_g u + u^3 - u,
$$
\noindent where $\Delta_g$ is the Laplacian operator of $g$. The Allen-Cahn Operator appears in mathematical physics to describe the process of phase separation in metal allows (c.f. \cite{AllenCahn}), and its interesting properties have already made it the object of various mathematical studies (c.f., for example, \cite{Cabre}, \cite{DeGiorgi} and \cite{Ilmanen}). In particular, the Allen-Cahn Operator is variational, arising as the Euler-Lagrange Equations (that is, the $L^2$-gradient) of the Ginzburg-Landau-Wilson Free Energy Functional:
$$
\Cal{E}_{\epsilon,g}(u) = \int_M\epsilon\|\nabla^g u\|^2 + \frac{1}{4}(u^2 - 1)^2\opdVol,
$$
\noindent and for this reason naturally lends itself to analysis by Morse theoretical techniques. In this note, we use Morse Homology, we show how the space of solutions to the Allen-Cahn Equation bifurcates as $\epsilon$ becomes small. Indeed, the Morse Homology yields a lower bound for the number of solutions, which increases discretely as $\epsilon^{-1}$ crosses points of the spectrum of $-\Delta_g$, tending to infinity as $\epsilon$ tends to zero.
\medskip
\noindent We denote by $\Cal{M}$ the space of smooth Riemannian metrics over $M$, which we furnish with the topology of $C^\infty$ convergence. For any smooth metric $g$, we define the {\bf solution space} $\Cal{Z}_g\subseteq]0,\infty[\times C^\infty(M)$ by:
$$
\Cal{Z}_g = \left\{(\epsilon, u)\ |\ \opAC_{\epsilon,g}(u) = 0\right\},
$$
\noindent so that $\Cal{Z}_g$ is the union of all solution sets for the metric $g$ and all parameters $\epsilon\in]0,\infty[$. We denote by $e_g:\Cal{Z}_g\rightarrow]0,\infty[$ the projection onto the second factor, and we will see presently (c.f. Proposition \procref{PropFirstPropernessResult}) that this is a proper map. For all $\epsilon\in]0,\infty[$, we denote:
$$
\Cal{Z}_{\epsilon,g} = e_g^{-1}(\left\{\epsilon\right\}) = \left\{u\ |\ \opAC_{\epsilon,g}(u) = 0\right\},
$$
\noindent so that $\Cal{Z}_{\epsilon,g}$ is the solution set for the metric $g$ and the parameter $\epsilon$. We aim to study the manner in which $\Cal{Z}_{\epsilon,g}$ bifurcates as $\epsilon$ tends to $0$, and the simplest way to do so is to study the geometry of $(\Cal{Z}_g,e_g)$. We first show that, upon perturbing the metric by an arbitrarily small amount, we may suppose that this geometry is relatively straightforward but for a countable set of singularities determined by the Laplacian of $g$. Indeed, we denote by $\opSpec(-\Delta_g)$ the set of all eigenvalues of $-\Delta_g$ (which, by convention, is non-negative), and we define the {\bf singular set}, $\opSing_g\subseteq]0,\infty[\times C^\infty(M)$ by:
$$
\opSing_g = \left\{(\epsilon,0)\ |\ \epsilon^{-1}\in\opSpec(-\Delta_g)\right\}.
$$
\noindent We recall that a subset $X$ of $\Cal{M}$ is said to be {\bf generic} (or equivalently, in the second category in the sense of Baire), whenever it contains a countable intersection of dense open sets. A given property is then said to hold for generic elements whenever it holds for all elements of some generic set. We recall that, by the Baire Category Theorem, any property that holds for generic elements of $\Cal{M}$ in particular holds over a dense subset of $\Cal{M}$. Using transversality techniques, we to show:
\proclaim{Theorem \nextprocno}
\noindent For generic $g\in\Cal{M}$, $\Cal{Z}_g\setminus\opSing_g$ is a smooth, $1$-dimensional submanifold of $]0,\infty[\times C^\infty(M)$. Moreover, if $\opDim(M)\geqslant 3$, then we may assume in addition that all critical points of $e_g$ are non-degenerate. In particular:
\medskip
\myitem{(1)} if $\epsilon^{-1}\notin\opSpec(-\Delta_g)$, then $\Cal{Z}_{\epsilon,g}$ is finite; and
\medskip
\myitem{(2)} there exists a discrete subset $X$ of the complement of $\opSpec(-\Delta_g)$ such that if $\epsilon^{-1}\notin X\munion\opSpec_g$, then $\Cal{Z}_{\epsilon,g}$ only consists of non-degenerate solutions of the Allen-Cahn Equation.
\endproclaim
\proclabel{ThmTheSolutionSpaceIsGenericallyNice}
\remark We recall that a solution to an elliptic partial differential equation is said to be non-degenerate whenever the linearisation of the operator about that solution is invertible.\qed
\medskip
\noindent When $u$ is a solution of the Allen-Cahn Equation with parameter $\epsilon$, we denote by $L\opAC_{\epsilon,g}(u)$ the linearisation of the Allen-Cahn Operator about $u$. When $u$ is non-degenerate, we denote by $\opIndex(u)$ its {\bf Morse Index}, which we recall is defined to be equal to the number of strictly negative eigenvalues of $L\opAC_{\epsilon,g}(u)$ counted with geometric multiplicity. Observe that the constant function $u=0$ is a solution of $\opAC_{\epsilon,g}$ for all $g$ and for all $\epsilon$. Moreover, as we will see presently (c.f. Proposition \procref{PropConstantSolutions}), $u=0$ is non-degenerate if and only if $\epsilon^{-1}\notin\opSpec(-\Delta_g)$, and the index of this solution is given by:
$$
\opIndex(0) = \#\left\{\lambda\in\opSpec(-\Delta_g)\ |\ \lambda <\epsilon^{-1}\right\}.
$$
\noindent Observe that $\opIndex(0)$ tends to infinity as $\epsilon$ tends to zero. Our second result now describes in terms of the number of solutions of a given Morse Index how $\Cal{Z}_{\epsilon,g}$ bifurcates as $\epsilon$ tends to $0$:
\proclaim{Theorem \nextprocno}
\noindent For generic $g\in\Cal{M}$, if $\epsilon^{-1}\notin\opSpec(-\Delta_g)$, then for all $0\leqslant k<\opIndex(0)$, there exist at least two non-degenerate solutions $u$ and $-u$ of the Allen-Cahn Equation such that $\opIndex(u)=\opIndex(-u)=k$.
\endproclaim
\proclabel{ThmBifurcation}
\remark In particular, we do not require Part $2$ of Theorem \procref{ThmTheSolutionSpaceIsGenericallyNice}. For countably many values of $\epsilon$, there may exist finitely many degenerate solutions. We simply ignore them.\qed
\medskip
\noindent Theorem \procref{ThmBifurcation} follows from Theorem \procref{ThmTheSolutionSpaceIsGenericallyNice} in a straightforward manner from standard Morse homological techniques. Indeed, for all $g$ and for all $\epsilon$, the constant functions $u=\pm 1$ are also solutions of the Allen-Cahn Equation, this time with Morse Index equal to $0$. Denoting $l=\opIndex(0)$, it follows that the chain groups $C_0$ and $C_l$ of the Morse-Complex of the Ginzburg-Landau-Wilson Free Energy Functional are at least $2$- and $1$- dimensional respectively. Since the underlying space (that is, $C^\infty(M)$) is contractible, the Morse-Homology, $H_k$, is non-trivial only for $k=0$. Finally, as the Allen-Cahn Operator is an odd operator, all the intermediate chain groups $C_1$, $C_2$, ...,$C_{l-1}$ have even dimension, and this fact, used together with the algebraic relations of Morse Homology allows us to deduce that they are non-trivial, thus proving the theorem.
\medskip
\noindent Theorem \procref{ThmTheSolutionSpaceIsGenericallyNice} is proven using the Sard-Smale Theorem, and this paper is therefore mostly devoted to obtaining the requisite surjectivity results. We draw the reader's attention to the fact that our usage of the Sard-Smale Theorem differs from standard approaches in one subtle but interesting respect. Indeed, whilst any application of the Sard-Smale Theorem is generally considered to require the separability of the function spaces used, we replace this condition by one that we call ``paraproperness'', which is to properness as paracompactness is to compactness. In Proposition \procref{PropPreservationOfParaproperness}, we make paraproperness into a useful concept by showing that it is preserved by restriction to both closed and open subsets, and in Theorem \procref{ThmSardSmale}, we reprove the Sard-Smale Theorem in this new context. This will be of particular use in the forthcoming paper \cite{SmiMH} where it makes possible the construction of a working Morse Homology theory in the H\"older space framework.
\medskip
\noindent This problem was recommended to the author by Frank Pacard. A large portion of this paper was written whilst the author was benefitting from a Marie Curie Postdoctoral Fellowship in the Centre de Recerca Matem\`atica, Barcelona, Spain. The author is also grateful to Joa Weber for encouragement and helpful suggestions to earlier drafts of this and the forthcoming paper.
\newhead{The Solution Space}
\newsubhead{Preliminaries and Compactness} For $\lambda\in[0,\infty]\setminus\Bbb{N}$, that is, for $\lambda=+\infty$, or for $\lambda=k+\alpha$, where $k\in\Bbb{N}$ and $\alpha\in]0,1[$, we denote by $C^\lambda:=C^\lambda(M)$ the space of $\lambda$-times H\"older differentiable functions over $M$, and when $\lambda<\infty$, we denote by $\|\cdot\|_\lambda$ the corresponding $C^\lambda$-H\"older norm. For $\mu\in[0,\infty]\setminus\Bbb{N}$, we likewise denote by $\Cal{M}^\mu$ the space of $C^\mu$-Riemannian metrics over $M$. It is well known that these spaces are non-separable, but as indicated in the introduction, this is of no consequence to us, and is satisfactorily treated by the concept of paraproperness (c.f. Section \subheadref{TheRegularSolutionSpace}, below).
\medskip
\noindent We consider a slightly more general problem than that discussed in the introduction. Let $f:\Bbb{R}\rightarrow\Bbb{R}$ be a smooth function such that $f$ is not linear over any interval, both $f$ and $f'$ have non-degenerate zeroes, and:
$$
\mlimsup_{t\rightarrow-\infty}f(t)<0,\qquad \mliminf_{t\rightarrow+\infty}f(t)>0.\eqnum{\nexteqnno}
$$
\eqnlabel{EqnFirstPropertiesOfPotentialFunction}%
\noindent As we will see presently (c.f. Proposition \procref{PropAPrioriBounds}, below) our theory only depends on the restriction of $f$ to the smallest interval containing all its zeroes. We therefore modify $f$ outside this interval, and replace \eqnref{EqnFirstPropertiesOfPotentialFunction} with the following technically more convenient property:
$$
\mlim_{t\rightarrow\pm\infty}f(t)/t\left|t\right| = +\infty.\eqnum{\nexteqnno}
$$
\noindent For $\mu>\lambda\in[0,\infty]\setminus\Bbb{N}$, we define the {\bf Allen-Cahn Operator}, $\opAC:]0,\infty[\times\Cal{M}^{\mu+1}\times C^{\lambda+2}\rightarrow C^\lambda$ by:
\eqnlabel{EqnPropertiesOfPotentialFunction}
$$
\opAC_{\epsilon,g}(u) := \opAC(\epsilon,g,u) = \epsilon\Delta_g u - f(u),\eqnum{\nexteqnno}
$$
\noindent where $\Delta_g$ is the Laplacian operator of $g$. Since $\opAC$ is constructed via a finite combination of multiplication, addition, differentiation and post-composition by smooth functions, it defines a smooth mapping between Banach manifolds. Importantly, the Allen-Cahn Operator arises as the $L^2$-gradient of the Ginzburg-Landau-Wilson Free Energy Functional:
$$
\Cal{E}_{\epsilon,g} = \int_M \epsilon\|\nabla^g u\|^2 + F(u)\opdVol,
$$
\noindent where $F$ is any primitive of $u$. In particular, solutions of the Allen-Cahn Equation are critical points of $\Cal{E}_{\epsilon,g}$.
\medskip
\noindent We define the {\bf solution space} $\Cal{Z}\subseteq]0,\infty[\times\Cal{M}^{\mu+1}\times C^{\lambda+2}$ by:
\eqnlabel{EqnAllenCahnEquation}
$$
\Cal{Z} = \opAC^{-1}(0).
$$
\noindent Let $\Pi:\Cal{Z}\rightarrow]0,\infty[\times\Cal{M}^{\mu+1}$ be the projection onto the first two factors and let $\Pi_g:\Cal{Z}\rightarrow\Cal{M}^{\mu+1}$ and $\Pi_u:\Cal{Z}\rightarrow C^{\lambda+2}$ be the projection onto the second and third factor, respectively. For all $(\epsilon,g)\in]0,\infty[\times\Cal{M}^{\mu+1}$, we define $\Cal{Z}_{\epsilon,g}\subseteq\Cal{Z}$, the {\bf solution space} for the data $(\epsilon,g)$ by:
$$
\Cal{Z}_{\epsilon,g} = \Pi^{-1}((\epsilon,g)).
$$
\noindent We study the bifurcations of $\Cal{Z}_{\epsilon,g}$ as $\epsilon$ varies. For this reason, we prefer to study all values of $\epsilon$ simultaneously and thus define $\Cal{Z}_g\subseteq\Cal{Z}$ by:
$$
\Cal{Z}_g = \Pi_g^{-1}(g).
$$
\noindent The main results of this paper follow from the differential topological properties of $\Cal{Z}$, $\Pi$ and $\Pi_g$, which we now proceed to study.
\medskip
\noindent We first review the analytic properties of the Allen-Cahn Operator. Elements of $\Cal{Z}$ have the following regularity properties:
\proclaim{Proposition \nextprocno}
\noindent Given $\mu>\lambda\in [0,\infty]\setminus\Bbb{N}$, if $(\epsilon,g,u)\in\Cal{Z}$, then $u\in C^{\mu+2}$.
\endproclaim
\proclabel{PropRegularityOfSolutions}
\proof Observe that $\epsilon\Delta_g$ is a second-order elliptic partial differential operator with coefficients in $C^{\mu}$. Thus, if $u$ lies in $C^{\mu+2(1-k)}$ for some positive integer $k$ with $\mu+2(1-k)>0$, then, since $f$ is smooth:
$$
\epsilon\Delta_g u = f\circ u \in C^{\mu+2(1-k)},
$$
\noindent and by elliptic regularity (c.f. \cite{GilbTrud}), $u\in C^{\mu+2(2-k)}$. Observe that since $u\in C^{\lambda+2}$, there exists $k$ such that $u\in C^{\mu+2(1-k)}$, and it follows by induction that $u\in C^{\mu+2}$, as desired.\qed
\medskip
\noindent In order to obtain a-priori estimates, we define $T_0>0$ by:
$$
T_0 = \msup\left\{ \left|t\right|\ |\ f(t)=0\right\}.
$$
\noindent It follows from \eqnref{EqnPropertiesOfPotentialFunction} that $T_0$ is finite. We have:
\proclaim{Proposition \nextprocno}
\noindent For all $(\epsilon,g,u)\in\Cal{Z}$:
$$
\|u\|_{L^\infty} \leqslant T_0.
$$
\endproclaim
\proclabel{PropAPrioriBounds}
\proof Suppose the contrary, that is, $\|u\|_{L^\infty}>T_0$. Since $M$ is compact, there exists $p\in M$ such that $\left|u(p)\right|=\|u\|_{L^\infty}$. If $u(p)\geqslant 0$, then $u(p)=\|u\|_{L^\infty}$, and since $p$ is a maximum of $u$, $(\Delta_gu)(p)\leqslant 0$, so that:
$$
f(\|u\|_{L^\infty}) = \epsilon(\Delta_g u)(p) \leqslant 0.
$$
\noindent On the other hand, if $u(p)<0$, then $u(p)=-\|u\|_{L^\infty}$ and $(\Delta_g u)(p)\geqslant 0$ so that:
$$
f(-\|u\|_{L^\infty}) = \epsilon(\Delta_g u)(p) \geqslant 0.
$$
\noindent In each case, this is absurd by definition of $T_0$ and Property \eqnref{EqnPropertiesOfPotentialFunction} of $f$, and the result follows.\qed
\proclaim{Proposition \nextprocno}
\noindent $\Pi$ defines a proper map from $\Cal{Z}$ into $]0,\infty[\times\Cal{M}^{\mu+1}$.
\endproclaim
\proclabel{PropFirstPropernessResult}
\proof Let $(\epsilon_m,g_m,u_m)_\minn$ be a sequence in $\Cal{Z}$ and suppose that $(\epsilon_m,g_m)_\minn$ converges to $(\epsilon_\infty,g_\infty)\in]0,\infty[\times\Cal{M}^{\mu+1}$, say. By Proposition \procref{PropRegularityOfSolutions}, $u_m\in C^{\mu+2}(M)$ for all $M$. By the Schauder estimates (c.f. \cite{GilbTrud}), there exists $B_1>0$ such that for all $m$:
$$\matrix
\|u_m\|_{\mu+2} \hfill& \leqslant B_1(\|u_m\|_{L^\infty} + \|\epsilon_m\Delta_{g_m}u_m\|_\mu)\hfill\cr
&= B_1(\|u_m\|_{L^\infty} + \|f\circ u_m\|_\mu)\hfill\cr
\endmatrix$$
\noindent By Proposition \procref{PropAPrioriBounds}, for all $m$, $u_m$ takes values in the compact set $[-T_0,T_0]$. Since $f$ is smooth, it follows from the chain rule and Gagliardo-Nirenberg-Moser type interpolation estimates (c.f. \cite{Taylor}) that there exists $B_2>0$ such that for all $m$:
$$
\|f\circ u_m\|_\mu \leqslant B_2(\|u_m\|_{L^\infty} + \|u_m\|_\mu).
$$
\noindent By standard interpolation inequalities (c.f. \cite{GilbTrud}), there exists $B_3\geqslant 0$ such that for all $m$:
$$
\|u_m\|_\mu \leqslant B_3\|u_m\|_{L^\infty} + \frac{1}{2B_1B_2}\|u_m\|_{\mu+2}.
$$
\noindent Combining these estimates yields, for all $m$:
$$\matrix
\|u_m\|_{\mu + 2}\hfill&\leqslant 2B_1(1+B_2)(1+B_3)\|u_m\|_{L^\infty}\hfill\cr
&\leqslant 2B_1(1+B_2)(1+B_3)T_0.\hfill\cr
\endmatrix$$
\noindent It now follows by the Arzela-Ascoli Theorem that there exists $u_\infty\in C^{\lambda+2}(M)$ towards which $(u_m)_\minn$ subconverges, and this completes the proof.\qed
\newsubhead{The Regular Solution Space and Paraproperness} For all $(\epsilon,g,u)\in\Cal{Z}$, we denote by $L\opAC$ the linearisation of $\opAC_{\epsilon,g}$ about $u$. By definition, $L\opAC=D_3\opAC(\epsilon,g,u)$, where $D_3\opAC$ denotes the partial derivative of $\opAC$ with respect to the third component. In particular, for all $(\epsilon,g,u)\in\Cal{Z}$ and for all $\varphi\in C^{\lambda+2}(M)$:
\subheadlabel{TheRegularSolutionSpace}
$$
L\opAC\varphi = \epsilon\Delta_g\varphi - f'(u)\varphi,\eqnum{\nexteqnno}
$$
\noindent so that $L\opAC$ is a self-adjoint second-order elliptic linear operator. In particular, it is Fredholm of index zero and by classical spectral theory, its spectrum is discrete, real and bounded above, and all of its eigenvalues have finite multiplicity. We say that the solution $u$ is {\bf non-degenerate} whenever $L\opAC$ is invertible, and we define the {\bf Morse Index} of $u$, which we denote by $\opIndex(u)$ by:
\eqnlabel{EqnLinearisedAllenCahnEquation}
$$
\opIndex(u) := \opIndex(L\opAC) = \sum_{\lambda\in\opSpec(L\opAC),\lambda>0}\opMult(\lambda),
$$
\noindent where, for all $\lambda\in\opSpec(L\opAC)$, $\opMult(\lambda)$ is its multiplicity.
\proclaim{Proposition \nextprocno}
\noindent For all $(\epsilon,g)\in]0,\infty[\times\Cal{M}^{\mu+1}$, the constant function $u=c$ is a solution to the Allen-Cahn Equation $\opAC_{\epsilon,g}(u)=0$ if and only if $f(c)=0$. Moreover, this solution is non-degenerate if and only if:
$$
\epsilon^{-1}f'(c) \notin \opSpec(\Delta_g),
$$
\noindent in which case its Morse Index is given by:
$$
\opIndex(c) = \sum_{\lambda\in\opSpec(\Delta_g),\lambda>\epsilon^{-1}f'(c)}\opMult(\lambda).
$$
\endproclaim
\proclabel{PropConstantSolutions}
\proof The first assertion follows immediately from \eqnref{EqnAllenCahnEquation}. By definition, $u$ is non-degenerate if and only if $0\notin\opSpec(L\opAC)$. By \eqnref{EqnLinearisedAllenCahnEquation}, this holds if and only if $\epsilon^{-1}f'(c)\notin\opSpec(L\opAC)$ and the second assertion follows. Finally, by \eqnref{EqnLinearisedAllenCahnEquation}, $\lambda>0$ is an eigenvalue of $L\opAC$ at $c$ if and only if $\mu:=\lambda+\epsilon^{-1}f(c)$ is an eigenvalue of $\Delta_g$. The third assertion then follows, and this completes the proof.\qed
\medskip
\noindent Since $f$ has non-degenerate zeroes, if $f(c)=0$, $f'(c)$ is either positive or negative. When $f'(c)$ is positive, it follows from Proposition \procref{PropConstantSolutions} that $u=c$ is always non-degenerate with Morse Index zero. On the other hand, if $f'(c)$ is negative, then $u=c$ is degenerate for countably many values of $c$ and its Morse Index tends to $+\infty$ as $\epsilon$ tends to $0$. It follows that zeroes of $f$ with negative derivative behave qualitatively differently from zeroes of $f$ with positive derivative. In fact, they yield singularities which are fundamental in the sense that they cannot be removed by perturbations of the metric, and this will be key to the bifurcation theory that follows. We therefore define the {\bf singular set}, $\opSing\subseteq]0,\infty[\times\Cal{M}^{\mu+1}\times\Cal{C}^{\lambda+2}$ by:
$$
\opSing = \left\{ (\epsilon,g,c)\ |\ f(c) = 0,\ \epsilon^{-1}f'(c)\in\opSpec(\Delta_g)\right\},
$$
\noindent and we define the {\bf regular solution space}, $\Cal{Z}^*\subseteq\Cal{Z}$ by:
$$
\Cal{Z}^* = \Cal{Z}\setminus\opSing.
$$
\noindent We now construct a countable exhaustion of $\Cal{Z}^*$ by closed sets. For all $g\in\Cal{M}^{\mu+1}$, we define:
$$
\opSing_g = \left\{ (\epsilon,c)\ |\ f(c) = 0,\ \epsilon^{-1}f'(c)\in\opSpec(\Delta_g)\right\}.
$$
\noindent By classical perturbation theory (c.f. \cite{Kato}), $\opSing_g$ varies continuously with $g$ in the Hausdorff sense. For all $m\in\Bbb{N}$, we define $\Cal{Z}_m\subseteq\Cal{Z}$ by:
$$
\Cal{Z}_m = \left\{(\epsilon,g,u)\in\Cal{Z}\ |\ 1/m\leqslant\epsilon\leqslant n,\ d((\epsilon,u),\opSing_g)\geqslant1/m\right\},
$$
\noindent and it follows from the continuous dependence of $\opSing_g$ on $g$ that $\Cal{Z}_m$ is closed. Moreover:
$$
\Cal{Z}^* = \munion_{m\in\Bbb{N}}\Cal{Z}_m.
$$
\noindent We now say that a continuous mapping $\Phi:X\rightarrow Y$ between two topological spaces is {\bf paraproper} whenever there exists a countable exhaustion $(X_m)_\minn$ of $X$ by closed sets such that for all $m$, the restriction of $\Phi$ to $X_m$ is proper. Paraproperness is made workable as a concept by the following restriction property:
\proclaim{Proposition \nextprocno}
\noindent Let $X$ and $Y$ be topological spaces and let $\Phi:X\rightarrow Y$ be paraproper.
\medskip
\myitem{(1)} if $K\subseteq X$ is closed, then the restriction of $\Phi$ to $K$ is paraproper; and
\medskip
\myitem{(2)} if $X$ is metrisable and if $\Omega\subseteq X$ is open, then the restriction of $\Phi$ to $\Omega$ is paraproper.
\endproclaim
\proclabel{PropPreservationOfParaproperness}
\proof Indeed, $(X_m\minter K)_\minn$ is a countable exhaustion of $K$ by closed sets and for all $m$, the restriction of $\Phi$ to $X_m\minter K$ is proper, which proves $(1)$. Now let $d$ be a distance function over $X$. For all $n\in\Bbb{N}$, we define $\Omega_n\subseteq X$ by:
$$
\Omega_n = \left\{ x\in X\ |\ d(x,\Omega^c)\geqslant (1/n)\right\}.
$$
\noindent For all $n$, $\Omega_n$ is closed, and since $\Omega$ is open:
$$
\Omega = \munion_{n\in\Bbb{N}}\Omega_n.
$$
\noindent $(\Omega_n\minter X_m)_{m,n\in\Bbb{N}}$ therefore constitutes a covering of $\Omega$ by closed sets. Moreover, for all $m,n\in\Bbb{N}$, since $\Omega_n\minter X_m$ is a closed subset of $X_m$, the restriction of $\Phi$ to this set is proper, and the restriction of $\Phi$ to $\Omega$ is therefore paraproper, which proves $(2)$.\qed
\medskip
\noindent In particular $\Pi_g$ defines a para-proper map from $\Cal{Z}$ into $\Cal{M}^{\mu+1}$:
\proclaim{Proposition \nextprocno}
\noindent For all $n$, $\Pi_g$ defines a proper map from $\Cal{Z}_n$ into $\Cal{M}^{\mu+1}$.
\endproclaim
\proclabel{PropPiIsParaproper}
\proof Let $(\epsilon_m,g_m,u_m)_\minn$ be a sequence in $\Cal{Z}_n$ and suppose that $(g_m)_\minn$ converges to $g_\infty\in\Cal{M}^{\mu+1}$. Since $\epsilon_m\in[1/n,n]$ for all $n$, by the Heine-Borel Theorem, we may suppose that there exists $\epsilon_\infty\in[1/n,n]$ towards which $(\epsilon_m)_\minn$ converges. By Proposition \procref{PropFirstPropernessResult}, there exists $u_\infty\in C^{\lambda+2}(M)$ towards which $(u_m)_\minn$ subconverges. Since $\Cal{Z}_n$ is closed, $(\epsilon_\infty,g_\infty,u_\infty)\in\Cal{Z}_n$, and the result follows.\qed
\medskip
\noindent Paraproperness now substitutes separability in our version of the Sard-Smale Theorem (c.f. \cite{Smale}):
\proclaim{Theorem \nextprocno, {\bf Sard-Smale}}
\noindent If $X$ and $Y$ are smooth Banach manifolds, and if $\Phi:X\rightarrow Y$ is a smooth, paraproper Fredholm map, then the set of regular values of $\Phi$ is generic in $Y$.
\endproclaim
\proclabel{ThmSardSmale}
\remark As Smale's result often mystifies, it is worth underlining the straightforward idea behind it. Using Fredholm Theory and the Implicit Function Theorem for Banach manifolds we reduce the problem to one of smooth maps between finite dimensional manifolds, and the result then follows by the classical Sard Theorem.\qed
\medskip
\proof Let $(X_n)_\ninn$ be a countable exhaustion of $X$ by closed sets such that for all $n$, the restriction of $\Phi$ to $X_n$ is proper. For all $n$, we denote the restriction of $\Phi$ to $X_n$ by $\Phi_n$, and we denote the set of regular values of $\Phi_n$ in $Y$ by $Y_n$. Since $\Phi_n$ is proper, and since surjectivity of Fredholm maps is an open property, $Y_n$ is open for all $n$.
\medskip
\noindent We now show that $Y_n$ is dense in $Y$. Indeed, choose $y\in Y$. Since we are only concerned with a neighbourhood of $y$ in $Y$, without loss of generality, we may suppose that $Y$ is a Banach space and that $y=0$. Define $\Psi:X\times Y\rightarrow Y$ by $\Psi(\tilde{x},\tilde{y})=\Phi(\tilde{x})+\tilde{y}$. Now choose $x\in\Phi_n^{-1}(0)$. Since $\Phi$ is Fredholm, $D\Phi(x)$ is closed and has finite dimensional cokernel, which we denote by $E_x$. In particular, the restriction of $D\Psi(x,0)$ to $T_xX\times E_x$ is surjective, and since surjectivity of Fredholm maps is an open property, there exists a neighbourhood $U_x$ of $x$ in $X_n$ such that the restriction of $D\Psi(\tilde{x},0)$ to $T_{\tilde{x}}X\times E_x$ is surjective for all $\tilde{x}\in U_x$. Since $\Phi_n^{-1}(0)$ is compact, it may be covered by finitely many such open sets, and there therefore exists a finite-dimensional subspace $E\subseteq Y$ such that the restriction of $D\Psi(\tilde{x},0)$ to $T_{\tilde{x}}X\times E$ is surjective for all $\tilde{x}\in\Phi_n^{-1}(y)$. We now consider the restriction of $\Psi$ to $X\times E$ and we denote $Z=\Psi^{-1}(0)$. By the Implicit Function Theorem for Banach manifolds, there exists a neighbourhood $\Omega$ of $\Phi_n^{-1}(0)\times\left\{0\right\}$ in $Z$ which is a smooth finite-dimensional submanifold of $X\times E$. Moreover, since $\Phi_n^{-1}(0)$ is compact, upon reducing $\Omega$ is necessary, we may suppose that this submanifold is separable. Let $\pi:\Omega\rightarrow E$ be the projection onto the first factor. Observe that if $\tilde{y}\in E$ is a regular value of $\pi$, then it is also a regular value of $\Phi_n$. However, by Sard's Theorem, regular values of $\pi$ are dense in $E$. It follows that $y=0$ is a concentration point of regular values of $\Phi_n$, and $Y_n$ is therefore a dense subset of $Y$ as asserted.
\medskip
\noindent Since the set of regular values of $\Phi$ coincides with $\minter_{n\in\Bbb{N}}Y_n$, it follows that this set is generic, which completes the proof.\qed
\newsubhead{The Regular Solution Space} In this section, we prove the following:
\proclaim{Proposition \nextprocno}
\noindent If $\opDim(M)\geqslant 2$, then for all $\mu>\lambda\in[0,\infty[\setminus\Bbb{N}$, $\Cal{Z}^*$ is a smooth Banach manifold modelled on $\Bbb{R}\times\Cal{M}^{\mu+1}$. Moreover, $\Pi_g$ defines a smooth, paraproper Fredholm map from $\Cal{Z}^*$ into $\Cal{M}^\mu$ of Fredholm index equal to $1$.
\endproclaim
\proclabel{PropRegularSolutionSpaceIsSmooth}
\noindent We prove this result using the Implicit Function Theorem for Banach manifolds. It is thus necessary to show that the derivative of $\opAC$ is surjective at every point of $\Cal{Z}^*$. We denote by $D_1\opAC$, $D_2\opAC$ and $D_3\opAC$ the partial derivatives of $\opAC$ with respect to the first, second and third components in $]0,\infty[\times\Cal{M}^{\mu+1}\times C^{\lambda+2}$ respectively. We are interested in particular in $D_2\opAC$. The tangent space of $\Cal{M}^{\mu+1}$ at any point canonically identifies with the space of $C^{\mu+1}$ sections of $\opSymm(TM)$. We denote this space by $\Gamma^{\mu+1}:=\Gamma^{\mu+1}(\opSymm(TM))$ and we refer to elements therin as {\bf first order perturbations} of the metric. We then identify $C^{\mu+1}$ with a subspace of $\Gamma^{\mu+1}$ by identifying every $C^{\mu+1}$ function $f$ with the $C^{\mu+1}$ section $fg$, and this induces the orthogonal splitting, $\Gamma^{\mu+1} = \Gamma_0^{\mu+1}\oplus C^{\mu+1}$, where $\Gamma_0^{\mu+1}:=\Gamma^{\mu+1}(\opSymm_0(TM))$ is the space of trace-free sections of $\opSymm(TM)$. The first order perturbations arising from sections of $C^{\mu+1}$ are precisely the conformal perturbations of the metric. However, it turns out that the useful perturbations for us are those whose trace vanishes. Indeed, for $g\in\Cal{M}^{\mu+1}$, and for any first order perturbation $A$ of $g$, denoting by $\delta_A\Delta_g$ the resulting first order perturbation of $\Delta_g$, we obtain:
\proclaim{Proposition \nextprocno}
\noindent If $\opTr(A)=0$, then, viewing $A$ as a section of $\opEnd(TM)$, for all $\varphi\in C^{\lambda+2}$:
$$
(\delta_A\Delta_g)\varphi = -\nabla\cdot(A\nabla\varphi),
$$
\noindent where $\nabla$ and $\nabla\cdot$ are the gradient and divergence operators of $g$ respectively.
\endproclaim
\proclabel{PropPerturbedLaplacian}
\proof We denote respectively by $\delta_A\Omega_g$ and $\delta_A\opHess_g$ the first order perturbations resulting from $A$ of the Levi-Civita covariant derivative and the Hessian operator of $g$. The Koszul formula yields:
$$
{(\delta_A\Omega_g)^k}_{;ij} = \frac{1}{2}({A^k}_{i;j} + {A^k}_{j;i} - {A_{ij;}}^k),
$$
\noindent where indices are raised and lowered with respect to $g$. Thus:
$$\matrix
&\delta_A\opHess_g(u)_{ij}\hfill&=-\frac{1}{2}({A^k}_{i;j} + {A^k}_{j;i} - {A_{ij;}}^k)u_{;k}\hfill\cr
\Rightarrow\hfill&\delta_A\Delta_g(u)\hfill&=-{A_j}^i{u_{;i}}^j - {A_j}^{i;j}u_{;i} + \frac{1}{2}\opTr(A)^{;k}u_{;k}\hfill\cr
& &=-({A_j}^iu_{;i})^{;j} + \frac{1}{2}\opTr(A)^{;k}u_{;k},\hfill\cr
\endmatrix$$
\noindent and since $\opTr(A)=0$, the result follows.\qed
\medskip
\noindent We recall the following straightforward result:
\proclaim{Proposition \nextprocno}
\noindent Let $X$ be a set consisting of at least $n$ distinct points. Let $E$ be an $n$-dimensional subset of the space of real-valued functions over $X$. Then there exist $n$ points $p_1,...,p_n\in X$ such that the mapping $\opEval:E\rightarrow\Bbb{R}^n$ given by:
$$
\opEval(f)_k = f(p_k),
$$
\noindent is a linear isomorphism.
\endproclaim
\proclabel{PropSamplingPoints}
\noindent This allows us to prove the required surjectivity result:
\proclaim{Proposition \nextprocno}
\noindent If $\opDim(M)\geqslant 2$, if $(u,g,\epsilon)\in\Cal{Z}$ and if $u$ is non-constant, then $D\opAC$ is surjective at $(u,g,\epsilon)$.
\endproclaim
\proclabel{PropSurjectivity}
\proof Since $D_3\opAC=L\opAC$ is elliptic, it has finite-dimensional cokernel, which we denote by $E$. Since $D_3\opAC$ is self-adjoint with respect to the $L^2$-inner-product of $g$, for all $\varphi\in E$:
$$
D_3\opAC(\varphi) = -\epsilon\Delta_g\varphi + f'(u)\varphi = 0.
$$
\noindent Since $u$ is non-constant, there exists $p\in M$ such that $\nabla u(p)\neq 0$. Let $\Omega$ be a neighbourhood of $p$ in $M$ diffeomorphic to the unit ball in Euclidean space over which $\nabla u$ doesn't vanish. By Aronszajn's unique continuation theorem (c.f. \cite{Aronszajn}), no non-trivial element of $E$ vanishes over $\Omega$. Furthermore, since $f'$ has non-degenerate zeroes, $f'(u)$ does not vanish identically over $\Omega$, and therefore no non-zero element of $E$ restricts to a constant map over this set. Thus, by Proposition \procref{PropSamplingPoints}, there exist $p_1,...,p_m\in\Omega\setminus\left\{p\right\}$ such that the mapping $\alpha:C^\lambda\rightarrow\Bbb{R}^m$ given by:
$$
\alpha(\varphi)_k = \varphi(p) - \varphi(p_k),
$$
\noindent restricts to a bijection on $E$.
\medskip
\noindent For any vector $\xi:=(\xi_0,...,\xi_m)$ of functions in $C_0^\infty(\Omega)$ we define $\alpha_\xi:C^\lambda\rightarrow\Bbb{R}^m$ by:
$$
\alpha_\xi(\varphi)_i = \int_M(\xi_0 - \xi_i)\varphi\opdVol_g.
$$
\noindent If $\xi_0-\xi_k$ is sufficiently close to $\delta_p-\delta_{p_k}$ in the weak sense for all $k$, where $\delta_p$ and $\delta_{p_k}$ are the Dirac delta functions supported at $p$ and $p_k$ respectively, then $\alpha_\xi$ is close to $\alpha$ and, in particular, is invertible. It follows that if $F$ is the linear span of $(\xi_0-\xi_k)_{1\leqslant k\leqslant m}$, then the $L^2$-inner-product restricts to a non-degenerate bilinear form over $E\times F$. In particular, $\opDim(F)=\opDim(E)$ and:
$$
F\minter\opIm(D_3\opAC) = F\minter E^\perp = \left\{0\right\}.
$$
\noindent $F$ is therefore complementary to $\opIm(D_3\opAC)$ in $C^\lambda$. That is:
$$
C^\lambda = F\oplus\opIm(D_3\opAC).
$$
\noindent However, we may suppose in addition that for all $1\leqslant k\leqslant m$:
$$
\int_M(\xi_0 -\xi_k)\opdVol_g=0
$$
\noindent It then follows from classical de-Rham cohomology theory that for all $k$ there exists a smooth vector field $X_k$ supported in $\Omega$ such that:
$$
\nabla\cdot X_k = \xi_0 - \xi_k.
$$
\noindent By Proposition \procref{PropRegularityOfSolutions}, $\nabla u$ is of class $C^{\mu+1}$, and thus, since it is non-vanishing over $\Omega$, there exists for all $k$ a $C^{\mu+1}$ field $A_k$ of symmetric matrices such that $A_k\nabla u = X_k$. In addition, since $M$ has dimension at least $2$, we may assume moreover that $\opTr(A_k)=0$ for all $k$, and it follows from Proposition \procref{PropPerturbedLaplacian} that:
$$
D_2\opAC\cdot A_k = -\epsilon\nabla\cdot(A_k\nabla u) = \epsilon\nabla\cdot X_k = \epsilon(\xi_0 - \xi_k).
$$
\noindent It follows that $F\subseteq\opIm(D_2\opAC)$ and so $C^\lambda\subseteq\opIm(D\opAC)$ and surjectivity follows.\qed
\medskip
\noindent Proposition \procref{PropRegularSolutionSpaceIsSmooth} follows readily:
\medskip
{\bf\noindent Proof of Proposition \procref{PropRegularSolutionSpaceIsSmooth}:\ }By Propositions \procref{PropConstantSolutions} and \procref{PropSurjectivity}, $D\opAC$ is surjective at every point of $\Cal{Z}^*$. Since $D_3\opAC$ is self-adjoint and elliptic, it is Fredholm of index zero, and it follows from the Implicit Function Theorem for Banach manifolds that $\Cal{Z}^*$ is a smooth Banach manifold modelled on $\Bbb{R}\times\Cal{M}^{\mu+1}$ and $\Pi_g$ is a smooth Fredholm map of Fredholm index equal to $1$. Finally, by Proposition \procref{PropPiIsParaproper}, $\Pi_g$ is para-proper, and this completes the proof.\qed
\medskip
\noindent Applying the Sard/Smale Theorem, now yields:
\proclaim{Proposition \nextprocno}
\noindent If $\opDim(M)\geqslant 2$, then for generic $g\in\Cal{M}^{\mu+1}$, $\Cal{Z}_g^*$ is a smooth, $1$-dimensional submanifold of $]0,\infty[\times C^{\lambda+2}$. Moreover, if we denote by $\epsilon_g:\Cal{Z}_g^*\rightarrow]0,\infty[$ the projection onto the first factor, then $\epsilon_g$ is proper.
\endproclaim
\proclabel{PropSmoothnessForGenericMetrics}
\remark Observe that paraproperness allows us to show that $\Cal{Z}_g^*$ is separable even though $\Cal{Z}^*$ isn't.\qed
\medskip
\proof By the Sard-Smale Theorem, the set of regular values of $\Pi_g$ is generic in $\Cal{M}^{\mu+1}$. Let $g\in\Cal{M}^{\mu+1}$ be a regular value of $\Pi_g$. By definition, $D\Pi_g(\epsilon,u,g)$ is surjective for all $(\epsilon,u)\in\Cal{Z}_g^*$. Since $\Pi_g$ is a smooth Fredholm map of Fredholm index equal to $1$, it follows from the Implicit Function Theorem for Banach manifolds that $\Cal{Z}_g^*$ is a (not necessarily separable) smooth, $1$-dimensional submanifold of $]0,\infty[\times C^{\lambda+2}$. By Proposition \procref{PropFirstPropernessResult}, $\epsilon_g:\Cal{Z}_g\rightarrow]0,\infty[$ is proper, and since $]0,\infty[$ has a compact exhaustion, so too does $\Cal{Z}_g$. In particular, $\Cal{Z}_g$ is separable, and therefore so too is $\Cal{Z}^*_g$, which completes the proof.\qed
\newsubhead{Non-Degeneracy of Critical Points} Let $\epsilon:]0,\infty[\times\Cal{M}^{\mu+1}\times C^{\lambda+2}\rightarrow]0,\infty[$ be the projection onto the first factor, and denote its restriction to $\Cal{Z}_g^*$ by $\epsilon_g$. We now aim to show that for generic $g\in\Cal{M}^{\mu+1}$, every critical point $\epsilon_g$ is non-degenerate. We first characterise those points where $d\epsilon_g$ vanishes:
\proclaim{Proposition \nextprocno}
\noindent At every point of $\Cal{Z}^*$:
$$
\opKer(D\epsilon)\minter\opKer(D\Pi_g)\minter T\Cal{Z}^* = \left\{0\right\}\times\left\{0\right\}\times\opKer(L\opAC).
$$
\endproclaim
\proclabel{PropFirstPropositionConcerningKernels}
\proof Indeed, by definition:
$$\matrix
T\Cal{Z}^*\hfill&=\opKer(D\opAC)\hfill\cr
&=\opKer(D_1\opAC\circ D\epsilon +  D_2\opAC\circ D\Pi_g + D_3\opAC\circ D\Pi_u),\hfill\cr
\endmatrix$$
\noindent where $D_1\opAC$, $D_2\opAC$ and $D_3\opAC$ represent the partial derivatives of $\opAC$ with respect to the first, second and third components respectively. Thus:
$$\matrix
\opKer(D\epsilon)\minter\opKer(D\Pi_g)\minter T\Cal{Z}^*
\hfill&=\opKer(D\epsilon)\minter\opKer(D\Pi_g)\minter\opKer(D_3\opAC\circ D\Pi_u)\hfill\cr
&=\left\{0\right\}\times\left\{0\right\}\times\opKer(L\opAC),\hfill\cr
\endmatrix$$
\noindent as desired.\qed
\proclaim{Proposition \nextprocno}
\noindent If $g$ is a regular value of $\Pi_g$, then at every point of $\Cal{Z}_g^*$:
$$
\opKer(d\epsilon_g) = \opKer(D\epsilon)\minter T\Cal{Z}^*_g = \left\{0\right\}\times\left\{0\right\}\times\opKer(L\opAC).
$$
\noindent In particular, $L\opAC$ has nullity at most $1$.
\endproclaim
\proclabel{PropSecondPropositionConcerningKernels}
\proof If $g$ is a regular value of $\Pi_g$, then:
$$
\opKer(D\Pi_g)\minter T\Cal{Z}^* = T\Cal{Z}^*_g,
$$
\noindent and the result now follows by Proposition \procref{PropFirstPropositionConcerningKernels}.\qed
\medskip
\noindent By Proposition \procref{PropSecondPropositionConcerningKernels}, if $g$ is a regular value of $\Pi_g$ and if $p\in\Cal{Z}_g$ is such that $d\epsilon_g=0$, then $\opKer(L\opAC)$ is $1$-dimensional. In particular, we may split $C^{\lambda+2}$ as the direct sum of $\opKer(L\opAC)$ and $\opKer(L\opAC)^\perp$ where $\opKer(L\opAC)^\perp$ is the orthogonal complement of $\opKer(L\opAC)$ in $C^{\lambda+2}$ with respect to the $L^2$-inner-product.
\proclaim{Proposition \nextprocno}
\noindent If $g$ is a regular value of $\Pi_g$, and if $p\in\Cal{Z}_g$ is such that $d\epsilon_g(p)=0$, then there is a neighbourhood $\Omega$ of $p$ in $\Cal{Z}^*$ which is a graph over $\Cal{M}^{\mu+1}\times\opKer(L\opAC)$.
\endproclaim
\proclabel{PropSolutionSpaceIsLocallyAGraph}
\proof Let $\pi:C^{\lambda+2}\rightarrow\opKer(L\opAC)$ be the orthogonal projection. Consider the restriction of the mapping $(\Pi_g,\pi\circ\Pi_u)$ to $\Cal{Z}^*$. Since $d\epsilon_g(p)=0$, bearing in mind Proposition \procref{PropSecondPropositionConcerningKernels}, at $p$:
$$
\opKer(D\Pi_g)\minter T\Cal{Z}^* = T\Cal{Z}_g^* = \opKer(d\epsilon_g)\minter T\Cal{Z}_g^* = \left\{0\right\}\times\left\{0\right\}\times\opKer(L\opAC).
$$
\noindent In particular, the restriction of $\pi\circ\Pi_u$ to $\opKer(D\Pi_g)$ is a linear isomorphism. The restriction of $(\Pi_g,\pi\circ\Pi_u)$ to $T\Cal{Z}^*$ is therefore also a linear isomorphism at $p$ and the result now follows by the Inverse Function Theorem for smooth maps between Banach manifolds.\qed
\medskip
\noindent Let $\Omega\subseteq\Cal{Z}^*$ be as in Proposition \procref{PropSolutionSpaceIsLocallyAGraph}. We construct a non-vanishing vector field, $X$ over $\Omega$ which is always tangent to $\Cal{Z}_g$ as follows. Choose $\varphi_0\in\opKer(L\opAC)$ such that $\|\varphi\|_{L^2}^2=1$. Let $X$ be the unique, smooth vector field over $\Omega$ which projects down to $\varphi$. There exist smooth functions $s:\Omega\rightarrow\Bbb{R}$ and $\varphi:\Omega\rightarrow\varphi_0+\opKer(L\opAC)^\perp$ such that, throughout $\Omega$:
$$
X = (s,0,\varphi).
$$
\noindent Trivially:
$$
D\Pi_g\cdot X = 0,
$$
\noindent so that $X$ is always tangent to $\Cal{Z}_g$, as desired.
\medskip
\noindent We now recall the following formula for the variation of a non-degenerate eigenvalue. Let $E\subseteq F\subseteq L^2(M)$ be Banach spaces and let $i:E\rightarrow F$ be a continuous embedding with dense image. It is normal to suppress $i$ and identify elements of $E$ with their image in $F$. Let $A\in\opLin(E,F)$ be a bounded, linear map. We recall that $E$ is said to be {\bf self-adjoint} if and only if for all $u,v\in E$:
$$
\langle u,A(v)\rangle = \langle A(u),v\rangle.
$$
\noindent The Implicit Function Theorem for Banach manifolds readily yields:
\proclaim{Proposition \nextprocno}
\noindent Let $X$, $E$ and $F$ be Banach spaces. Let $A:X\rightarrow\opLin(E,F)$ be a smooth mapping such that for all $x\in X$, $A_x:=A(x)$ is self-adjoint and Fredholm of index zero. Suppose that $\opNull(A_0)=1$ and let $\varphi_0$ be a non-zero element of $\opKer(A_0)$. Then there exists a neighbourhood $U$ of $0$ in $X$ and smooth maps $\lambda:X\rightarrow\Bbb{R}$ and $\varphi:X\rightarrow\varphi+\opKer(A_0)^\perp$ such that $\lambda(0)=0$, $\varphi(0)=\varphi_0$ and for all $x\in X$:
$$
A(x)\varphi(x) = \lambda(x)\varphi(x).
$$
\noindent Moreover, for any tangent vector $\xi$ to $X$ at $0$:
$$
d\lambda(\xi) = \langle DA_0(\xi)\varphi_0,\varphi_0\rangle.
$$
\endproclaim
\proclabel{PropVariationOfNonDegenerateEigenvalue}
\noindent By Proposition \procref{PropVariationOfNonDegenerateEigenvalue}, upon reducing $\Omega$ if necessary, there exist smooth functions $\lambda:\Omega\rightarrow\Bbb{R}$ and $\tilde{\varphi}:\Omega\rightarrow\varphi_0+\opKer(L\opAC(p_0))^\perp$ such that $\lambda(p_0)=0$, $\tilde{\varphi}(p_0)=\varphi_0$ and throughout $\Omega$:
$$
L\opAC\tilde{\varphi}=\lambda\tilde{\varphi}.
$$
\noindent The role played by $\lambda$ is revealed by the following quantitative analogue of Proposition \procref{PropSecondPropositionConcerningKernels}:
\proclaim{Proposition \nextprocno}
\noindent Let $g\in\Cal{M}^{\mu+1}$ be a regular value of $\Pi_g$. If $p\in\Cal{Z}_g$ is such that $d\epsilon_g(p)=0$, then:
$$
\langle D_1\opAC\cdot D^2\epsilon_g(X_p,X_p),\varphi_0\rangle = -d\lambda(X_p).
$$
\noindent In particular, if $p$ is a non-degenerate zero of $\lambda$, then it is also a non-degenerate zero of $d\epsilon_g$.
\endproclaim
\proclabel{PropRelationBetweenLambdaAndEpsilon}
\proof By definition, $\opAC$ vanishes over $\Cal{Z}^*$ and so:
$$
D_1\opAC\circ D\epsilon(X) = -D_2\opAC\circ D\Pi_g(X) - D_3\opAC\circ D\Pi_u(X) = - L\opAC\circ D\Pi_u(X).
$$
\noindent Since $D\epsilon(X_p)=0$, differentiating a second time yields:
$$
D_1\opAC\circ D^2\epsilon(X_p,X_p) = -D_{X_p}L\opAC\circ D\Pi_u(X_p) - L\opAC\circ D^2\Pi_u(X_p,X_p).
$$
\noindent Observe that $D^2\Pi_u(X_p,X_p)$ takes values in $\opKer(L\opAC)^\perp$. Moreover, $\varphi_0\in\opKer(L\opAC)$, and since $L\opAC$ preserves both $\opKer(L\opAC)$ and $\opKer(L\opAC)^\perp$, taking the inner-product with $\varphi_0$ yields:
$$\matrix
\langle D_1\opAC\circ D^2\epsilon(X_p,X_p),\varphi_0\rangle\hfill
&=-\langle D_{X_p}L\opAC\circ D\Pi_u(X_p),\varphi_0\rangle\hfill\cr
&=-\langle D_{X_p}L\opAC\varphi_0,\varphi_0\rangle.\hfill\cr
\endmatrix$$
\noindent Thus, by Proposition \procref{PropVariationOfNonDegenerateEigenvalue}:
$$
\langle D_1\opAC\circ D^2\epsilon_g(X_p,X_p),\varphi_0\rangle = -d\lambda(X_p),
$$
\noindent as desired.\qed
\medskip
\noindent The above discussion is most usefully summarised as follows:
\proclaim{Proposition \nextprocno}
\noindent There exists an open subset $\Omega\subseteq\Cal{Z}^*$ and a smooth function $\lambda:\Omega\rightarrow\Bbb{R}$ with the following properties:
\medskip
\myitem{(1)} if $p\in\Cal{Z}^*$ is such that $D\Pi_g(p)$ is surjective and $d\epsilon_g(p)=0$, then $p\in\Omega$ and $\lambda=0$; and
\medskip
\myitem{(2)} for all $p\in\Omega$, $\lambda(p)$ is an eigenvalue of $L\opAC(p)$.
\endproclaim
\proclabel{PropLambdaDefinedGlobally}
\proof Let $p\in\Cal{Z}^*$ be such that $D\Pi_g(p)$ is surjective and $D\epsilon_g(p)=0$. Let $\Omega_p$ and $\lambda:\Omega_p\rightarrow\Bbb{R}$ be as in the preceeding discussion. Upon reducing $\Omega_p$ if necessary, we may assume that $\lambda$ is the eigenvalue of $L\opAC(p)$ with least absolute value. It follows then that $\lambda$ is uniquely defined, and taking the union over all such $\Omega_p$ yields the desired open set and smooth function.\qed
\proclaim{Proposition \nextprocno}
\noindent Suppose that $\opDim(M)\geqslant 3$. Choose $p\in\Omega$ and let $\varphi\in C^{\lambda+2}$ be an element of $\opKer(L\opAC(p))$. If there exists a point $x\in M$ such that $du(x)$ and $d\varphi(x)$ are both non-vanishing and non-colinear, then $d\lambda$ is non-zero at $p$.
\endproclaim
\proclabel{PropDerivativeOfLambdaIsNonVanishing}
\proof Let $A$ be a trace-free first order perturbation of $g$ and let $\delta_A\opdVol_g$, $\delta_AL\opAC$ and $\delta_A\lambda$ denote the resulting first order perturbations of $\opdVol_g$, $L\opAC$ and $\lambda$ respectively. Then:
$$
\delta_A\opdVol_g = \opTr(A)\opdVol_g = 0.
$$
\noindent Thus, by Proposition \procref{PropVariationOfNonDegenerateEigenvalue}:
$$
\delta_A\lambda = \int\varphi(\delta_AL\opAC)\varphi\opdVol_g.
$$
\noindent However, by Proposition \procref{PropPerturbedLaplacian}:
$$
\delta_AL\opAC\varphi = -\epsilon\nabla\cdot(A\nabla\varphi),
$$
\noindent and so:
$$\matrix
\delta_A\lambda \hfill&= -\epsilon\int\varphi\nabla\cdot(A\nabla\varphi)\opdVol_g\hfill\cr
&=\epsilon\int\langle A,\nabla\varphi\otimes\nabla\varphi\rangle\opdVol_g.\hfill\cr
\endmatrix$$
\noindent Let $p\in M$ be such that $du(p)$ and $d\varphi(p)$ are non-vanishing and non-colinear. Since $M$ is $3$-dimensional, there exists a first-order perturbation $A$ of $g$, supported near $p$ such that:
$$
\opTr(A) = 0,\qquad (\delta_A\Delta_g)u=0,\qquad \delta_A\lambda\neq 0.
$$
\noindent It follows from the first two relations that the vector $(0,A,0)$ is tangent to $\Cal{Z}^*$, and it follows from the third relation that $d\lambda(0,A,0)\neq 0$, which completes the proof.\qed
\medskip
\noindent Applying the Sard/Smale Theorem, we now obtain:
\proclaim{Proposition \nextprocno}
\noindent If $\opDim(M)\geqslant 3$, then for generic $g\in\Cal{M}^{\mu+1}$ and for $p\in\Cal{Z}_g$, if $d\epsilon_g(p)=0$, then either:
\medskip
\myitem{(1)} $D^2\epsilon_g(p)\neq 0$; or
\medskip
\myitem{(2)} if $\varphi\in\opKer(L\opAC(p))$, then for all $Y\in TM$, if $du(Y)=0$, then $d\varphi(Y)=0$.
\endproclaim
\proclabel{PropAlmostNonDegenerateEpsilonForGenericMetrics}
\proof Let $X\subseteq\Omega$ be the set of all points $p$ such that $d\epsilon_g(p)=0$ and $(2)$ is satisfied. Observe that $(2)$ implies that $du$ and $d\varphi$ are everywhere colinear. Since this is a closed condition, it follows that $X$ is a closed subset of $\Omega$, and $\tilde{\Omega}:=\Omega\setminus X$ is therefore open. Let $Y\subseteq\tilde{\Omega}$ be the set of all points where $\lambda$ vanishes. By Proposition \procref{PropDerivativeOfLambdaIsNonVanishing} and the Implicit Function Theorem for Banach manifolds, $Y$ is a smooth codimension-$1$ Banach submanifold of $\tilde{\Omega}$. Observe that the restriction of $\Pi_g$ to $Y$ is a smooth Fredholm map of Fredholm index $0$. Moreover, by Proposition \procref{PropPreservationOfParaproperness}, this restriction is paraproper. It therefore follows from Theorem \procref{ThmSardSmale} that for generic $g\in\Cal{M}^{\mu+1}$, $g$ is a regular value of this restriction. Moreover, since the intersection of two generic sets is also generic, we may assume that $g$ is also a regular value of $\Pi_g$. For such a $g$, $\Cal{Z}_g$ is a smooth $1$-dimensional manifold and the restriction of $\lambda$ to $\Cal{Z}_g$ has non-degenerate zeroes at all points where $(2)$ is satisfied, and the result now follows by Proposition \procref{PropRelationBetweenLambdaAndEpsilon}.\qed
\newsubhead{The Degenerate Case} We now eliminate Case $(2)$ of Proposition \procref{PropAlmostNonDegenerateEpsilonForGenericMetrics}. We begin by characterising its geometry:
\proclaim{Proposition \nextprocno}
\noindent Let $u,\varphi\in C^{\lambda+2}(M)$ be such that $u$ is non-constant, $\varphi$ is non-zero, and:
$$
\epsilon\Delta_g u = f(u),\qquad \epsilon\Delta_g\varphi = f'(u)\varphi.
$$
\noindent If for all vectors $X\in TM$ such that $du(X)=0$ we have $d\varphi(X)=0$, then $\|du\|_g$ is constant over each connected component of every level set of $u$.
\endproclaim
\proclabel{PropGeometricPropertiesOfDegenerateCase}
\remark In fact, we can prove more: the complement of the vanishing set of $du$ in $M$ is foliated by compact hypersurfaces of constant mean curvature. This property is interesting, as it is independent of the parameter $\epsilon$.\qed
\medskip
\proof Observe that if $u$ is constant over any non-trivial neighbourhood, then it is equal to a zero of $f$, $c$ say over this neighbourhood. Since $u=c$ is also a solution of $\opAC(u)=0$, it follows from Aronszajn's unique continuation theorem (c.f. \cite{Aronszajn}) that $f=c$ over the whole of $M$, which is absurd, and it follows that $du$ is almost everywhere non-vanishing.
\medskip
\noindent Now choose $p\in M$ such that $du(p)\neq 0$. Let $\Omega$ be a neighbourhood of $p$ over which $du$ does not vanish. Observe that the image of the restriction of $u$ to $\Omega$ is an open interval, $I$, say. Moreover, by Aronszajn's unique continuation theorem again, the restriction of $\varphi$ to $\Omega$ is non-zero. Let $\Cal{F}$ denote the foliation of $\Omega$ by level hypersurfaces of $u$. By hypothesis, $\varphi$ is constant over each leaf of $\Cal{F}$. Thus, upon reducing $\Omega$ if necessary, there exists a non-zero $C^{\lambda+2}$-function $\Phi:I\rightarrow\Bbb{R}$ such that, over $\Omega$, $\varphi = \Phi(u)$. Taking the Laplacian of both sides of this relation yields:
$$
f'(u)\Phi(u) = \epsilon\Phi''(u)\|du\|_g^2 + f(u)\Phi'(u).
$$
\noindent We claim that $\Phi''$ is almost everywhere non-vanishing. Indeed, otherwise, upon reducing $\Omega$ further if necessary, we may suppose that $\Phi$ is linear and that $f'\Phi - f\Phi'=0$. The restriction of $f$ to $I$ is therefore also linear, which is absurd by the hypothesis on $f$, and $\Phi''$ is therefore almost everywhere non-vanishing, as asserted. However, whenever $\Phi''(u)\neq 0$, we have:
$$
\|du\|_g^2 = \frac{1}{\epsilon\Phi''(u)}(f'(u)\Phi(u) - f(u)\Phi'(u)),
$$
\noindent from which it follows that $\|du\|_g^2$ is constant over every leaf of $\Cal{F}$ where $\Phi''(u)$ does not vanish. Since the set of all such leaves is dense, it follows that $\|du\|_g$ is constant over every leaf of $\Cal{F}$.
\medskip
\noindent Choose $t\in\Bbb{R}$ and denote $X=u^{-1}(t)$. Let $X_0,X_1\subseteq X$ be respectively the subset of $X$ consisting of those points where $du$ vanishes, and the subset of $X$ consisting of those points where it does not vanish. Trivially, $\|du\|_g$ is constant over $X_0$. Observe that $X_1$ is a submanifold of $M$. Moreover, by the above discussion, $\|du\|_g$ is constant over every connected component of $X_1$. Every connected component of $X_1$ is therefore a closed submanifold, and, in particular, is disjoint from $X_0$. It follows that if $X'$ is a connected component of $X$, then $X'$ is either contained wholly in $X_0$ or wholly in $X_1$. In either case, $\|du\|_g$ is constant over $X'$, and this completes the proof.\qed
\medskip
\noindent The following refinement of Proposition \procref{PropGeometricPropertiesOfDegenerateCase} is easier to work with:
\proclaim{Proposition \nextprocno}
\noindent Under the same hypotheses as Proposition \procref{PropGeometricPropertiesOfDegenerateCase}, if $X_p\in TM$ is such that $du(X_p)=0$, then:
$$
\opHess^g(u)(\nabla^g u,X_p) = 0.
$$
\endproclaim
\proclabel{PropSecondVersionGeometricPropertiesOfDegenerateCase}
\proof If $du(p)=0$, then $\nabla^g u(p)=0$, and the result follows trivially. Otherwise, $du(p)\neq 0$, and, by Proposition \procref{PropGeometricPropertiesOfDegenerateCase}, $\|\nabla^g u\|=\|du\|_g$ is constant over the level hypersurface of $u$ passing through $p$. Since $du(X_p)=0$, $X_p$ is tangent to $S$, and so:
$$
\opHess(u)(\nabla^g u, X_p) = \langle \nabla_{X_p}\nabla^g u,\nabla^g u\rangle = \frac{1}{2}X_p\|\nabla^g u\|^2 = 0,
$$
\noindent as desired.\qed
\proclaim{Proposition \nextprocno}
\noindent Suppose that $\opDim(M)\geqslant 2$ and let $u\in C^{\lambda+2}$ be a non-constant function such that $\epsilon\Delta_g u = f(u)$. Choose $p\in M$ such that $\nabla^g u(p)\neq 0$ and $Y\in T_pM$ such that $du(Y)=0$. There exists a $C^{\mu+1}$ first order perturbation $A$ of the metric supported in an arbitrarily small neighbourhood of $p$ such that such that:
\medskip
\myitem{(1)} $A(p)=0$;
\medskip
\myitem{(2)} $(\delta_A\Delta_g)u=0$; and
\medskip
\myitem{(3)} $\delta_A\opHess(u)(\nabla u,Y)(p)\neq 0$.
\endproclaim
\proclabel{PropRemovingGeodesicProperty}
\proof Let $\Omega$ be a neighbourhood of $p$ diffeomorphic to the unit ball. Let $X$ be a smooth divergence-free vector field supported in $\Omega$ such that $X(p)=0$. Since $M$ is at least two dimensional, and since $\nabla u$ does not vanish over $\Omega$, there exists a $C^{\mu+2}$ section $A$ of $\opSymm(TM)$ supported in $\Omega$ such that $A\cdot\nabla u=X$ and $\opTr(A)=0$. In particular, we may suppose that $A(p)=0$. By Proposition \procref{PropPerturbedLaplacian}, for any such $A$:
$$
(\delta_A\Delta_g)u=0.
$$
\noindent As in the proof of Proposition \procref{PropPerturbedLaplacian}, the first order perturbation of the Hessian of $u$ is given by:
$$
\delta_A\opHess_g(u)(X,Y) =\frac{1}{2}({A_{ij;}}^k - {A^k}_{i;j} - {A^k}_{j;i})u_kX^iY^j.
$$
\noindent Thus, bearing in mind that $(\delta_A\nabla^g)u(p)=X(p)=0$:
$$\matrix
\delta_A\opHess_g(u)(\nabla^g u, Y)\hfill&=\frac{1}{2}({A_{ij;}}^k - {A^k}_{i;j} - {A^k}_{j;i})u_ku^iY^j\hfill\cr
&=-\frac{1}{2}A_{ki;j}u^ku^iY^j\hfill\cr
&=-\frac{1}{2}Y\langle A\nabla^g u, \nabla^g u\rangle + \langle A\nabla^g_Y\nabla^g u,\nabla^g u\rangle\hfill\cr
&=-\frac{1}{2}Y\langle A\nabla^g u,\nabla^g u\rangle + \opHess_g(u)(Y,X)\hfill\cr
&=-\frac{1}{2}Y\langle X,\nabla^g u\rangle\hfill\cr
&=-\frac{1}{2}Ydu(X).\hfill\cr
\endmatrix$$
\noindent Since $X$ is divergence free and compactly supported in $\Omega$, it follows from classical de-Rham cohomology theory that there exists a $2$-form $Z$ supported in $\Omega$ such that:
$$
X = \nabla^g\cdot Z.
$$
\noindent We choose exponential coordinates about $p$, and write $Z$ as:
$$
Z = \sum_{i<j}Z^{ij}\partial_i\wedge\partial_j,
$$
\noindent so that:
$$
X^k = \sum_{i>k}\partial_iZ^{ik} - \sum_{i<k}\partial_iZ^{ki}.
$$
\noindent We choose the basis at $p$ such that $\nabla^g u$ and $Y$ are colinear with $\partial_1$ and $\partial_2$ respectively. Thus, at the origin, bearing in mind that $X(p)=0$:
$$
-\frac{1}{2}Ydu(X) = \frac{1}{2}\|\nabla^g u\|_g\|Y\|_g\sum_{i>1}\partial_2\partial_iZ^{1i}.
$$
\noindent We choose $Z$ such that $\partial_iZ^{jk}=0$ for all $i$, $j$ and $k$, $\partial_2\partial_1 Z^{11}=1$ and $\partial_2\partial_i Z^{1i}=0$ for all $i>1$. Then, if $X=\nabla^g\cdot Z$:
$$
X(p)=0,\qquad -\frac{1}{2} Ydu(X) = \frac{1}{2}\|\nabla^g u\|_g\|Y\|_g.
$$
\noindent $X=\nabla^g\cdot Z$ is the desired vector field, and this completes the proof.\qed
\proclaim{Proposition \nextprocno}
\noindent If $\opDim(M)\geqslant 2$, then for generic $g\in\Cal{M}^{\mu+1}$, if $(\epsilon,g,u)\in\Cal{Z}_g$, if $u$ is non-constant and if $\varphi\in\opKer(LAC(u))$ is non-zero, then there exists a point $p\in M$ such that $du(p)$ and $d\varphi(p)$ are both non-zero and non-colinear.
\endproclaim
\proclabel{PropDegenerateCaseDoesNotExistForGenericMetrics}
\proof Let $X_p$ and $Y_q$ be unit vectors over distinct points of $M$. Let $\Omega:=\Omega(X_p,Y_q)\subseteq\Cal{Z}^*$ be the open set of all $(\epsilon,g,u)$ such that $\nabla^g u(p)$ and $\nabla^g u(q)$ are both non-zero and non-colinear with $X_p$ and $Y_q$ respectively. We define the functions $\Phi_p,\Phi_q:\Omega\rightarrow\Bbb{R}$ by:
$$
\Phi_p(\epsilon,g,u) = \opHess_g(u)(X_p^\perp,\nabla^g u(p)),\qquad \Phi_q(\epsilon,g,u) = \opHess_g(u)(Y_q^\perp,\nabla^g u(q)),
$$
\noindent where $X_p^\perp$ and $Y_q^\perp$ are the orthogonal projections of $X_p$ and $Y_q$ respectively onto the normal hyperplanes to $\nabla^g u(p)$ and $\nabla^g u(q)$ respectively. Observe that both $\Phi_p$ and $\Phi_q$ define smooth functions over $\Omega$. Moreover, it follows from Proposition \procref{PropRemovingGeodesicProperty} that $D(\Phi_p,\Phi_q)$ is surjective at every point of $\Omega$. Thus, if $Z:=Z(X_p,Y_q)$ is the zero set of this functional then it is a smooth, codimension $2$ submanifold of $\Omega$. In particular, the restriction of $\Pi_g$ to $Z$ is a smooth Fredholm map of index $-1$. Thus, if $g\in\Cal{M}^{-1}$ is a regular value of the restriction of $\Pi_g$ to $Z$, then $\Pi_g^{-1}(g)\minter Z$ is a smooth submanifold of $Z$ of dimension equal to $-1$, that is, it is empty. However, by Proposition \procref{PropPreservationOfParaproperness}, the restriction of $\Pi_g$ to $\Omega$, and therefore also to $Z$, is para-proper, and it follows by Theorem \procref{ThmSardSmale} that the set of regular values of this restriction is generic in $\Cal{M}^{\mu+1}$.
\medskip
\noindent Let $\Cal{X}\subseteq (UM\times UM)\setminus\pi^{-1}(\opDiag)$ be a countable dense family of pairs $(X_p,Y_q)$ of unit vectors above distinct points of $M$. Since the intersection of a countable family of generic sets is generic, it follows that for generic $g\in\Cal{M}^{\mu+1}$, and for all $(X_p,Y_q)\in\Cal{X}$, $\Pi_g^{-1}(g)\minter Z(X_p,Y_q)$ is empty. For such a $g$, choose $(\epsilon,g,u)\in\Cal{Z}_g$. Let $\tilde{p},\tilde{q}\in M$ be distinct points such that both $du(\tilde{p})$ and $du(\tilde{q})$ are non-zero, and let $\tilde{X}_{\tilde{p}}$ and $\tilde{Y}_{\tilde{q}}$ be unit vectors in $UM$ normal to $\nabla^gu(\tilde{p})$ and $\nabla^gu(\tilde{q})$ respectively. Since $\Cal{X}$ is dense, there exists a pair $(X_p,Y_q)\in\Cal{X}$ such that $du(p)$ and $du(q)$ are non-zero and $X_p$ and $Y_q$ are non-colinear with $\nabla^g u(p)$ and $\nabla^g u(q)$ respectively. However, by definition of $g$, $(\epsilon,g,u)\notin Z(X_p,Y_q)$, from which it follows that one of $\Phi_p(\epsilon,g,u)$ and $\Phi_q(\epsilon,g,u)$ is non-zero. In other words, without loss of generality:
$$
\opHess_g(u)(X_p^\perp,\nabla^g u(p)) \neq 0,
$$
\noindent and it now follows from Proposition \procref{PropSecondVersionGeometricPropertiesOfDegenerateCase} that there exists at least one point in $M$ where $du$ and $d\varphi$ are non-zero and non-colinear, as desired.\qed
\medskip
\noindent Combining these relations, we obtain Theorem \procref{ThmTheSolutionSpaceIsGenericallyNice}:
\medskip
{\bf\noindent Proof of Theorem \procref{ThmTheSolutionSpaceIsGenericallyNice}:\ }Since the intersection of finitely many generic sets is generic, this follows from Propositions \procref{PropSmoothnessForGenericMetrics}, \procref{PropAlmostNonDegenerateEpsilonForGenericMetrics} and \procref{PropDegenerateCaseDoesNotExistForGenericMetrics}.\qed
\newsubhead{The Solution Space at Infinity} Now fix $g\in\Cal{M}^{\mu+1}$. We show that for $\epsilon$ sufficiently large, the only elements of $\Cal{Z}_{\epsilon,g}$ are the constant solutions. We recall that for all $g$, $\opKer(\Delta_g)^\perp$ coincides with the space of functions whose integral with respect to the volume form of $g$ vanishes.
\proclaim{Proposition \nextprocno}
\noindent Let $c\in\Bbb{R}$ be such that $f(c)=0$. There exist $B>0$ and $\delta>0$ such that if $\epsilon>B$, if  $v\in C^{\lambda+2}$ and $t\in\Bbb{R}$ are such that:
$$
\int_Mv\opdVol_g=0,\qquad \|v\|_{\lambda+2} < \delta\epsilon^{-1},\qquad \left|t\right|<\delta,
$$
\noindent and if $\opAC(\epsilon,g,c + v +t)=0$, then $(v,t)=(0,0)$.
\endproclaim
\proclabel{PropPreliminaryToCondensationOfSolutions}
\proof Define $\Cal{F}:\opKer(\Delta_g)^\perp\times\Bbb{R}^2\rightarrow C^\lambda$ by:
$$
\Cal{F}(v,t,\eta) = \Delta_g v - f(c + \eta v + t).
$$
\noindent Observe that $\Cal{F}$ is a smooth function between Banach manifolds. Moreover, if we denote by $D_1\Cal{F}$ and $D_2\Cal{F}$ its partial derivatives with respect to the first and second factors respectively, then since $f'(c)\neq 0$, $D_1\Cal{F}+D_2\Cal{F}$ is surjective at $(0,0,0)$. It follows from the Implicit Function Theorem for Banach manifolds that there exists $b>0$ and a neighbourhood $W$ of $(0,0)$ in $\opKer(\Delta_g)^\perp\times\Bbb{R}$ such that if $\eta<b$ then there exists a unique point $(v_\eta,t_\eta)\in W$ such that $\Cal{F}(v_\eta,t_\eta,\eta)=0$. Since, in particular, $\Cal{F}(0,0,\eta)=0$ for all $\eta$, it follows that if $(v,t)\in W$ is such that $\Cal{F}(v,t,\eta)=0$, then $(v,t)=(0,0)$. Let $B=1/b$ and let $\delta>0$ be such that:
$$
\left\{(v,t)\ |\ \|v\|_{\lambda+2}<\delta,\ \left|t\right|<\delta\right\} \subseteq W.
$$
\noindent We claim that $B$ and $\delta$ have the desired properties. Indeed, let $\epsilon>B$, $v\in C^{\lambda+2}$ and $t\in\Bbb{R}$ be such that $v\in\opKer(\Delta_g)^\perp$, $\|v\|_{\lambda+2}<\delta\epsilon^{-1}$, $\left|t\right|<\delta$ and $\opAC(\epsilon,g,c+v+t)=0$. Then, denoting $\eta=1/\epsilon$:
$$
\Cal{F}(\epsilon v,t,\eta) = \Delta_g(\epsilon v) - f(c + \eta(\epsilon v) + t) = 0.
$$
\noindent Since $\|\epsilon v\|_{\lambda+2},\left|t\right|<\delta$, it follows from the preceeding discussion that $(v,t)=(0,0)$, as desired.\qed
\proclaim{Proposition \nextprocno}
\noindent There exists $B>0$ such that if $\epsilon>B$, then $\Cal{Z}_{\epsilon,g}$ only consists of constant solutions.
\endproclaim
\proclabel{PropCondensationOfSolutions}
\proof Suppose the contrary. There exists a sequence $(u_n,t_n,\epsilon_n)_\ninn\in\opKer(\Delta_g)^\perp\times\Bbb{R}^2$ such that $(\epsilon_n)_\ninn$ tends to $+\infty$, $u_n$ is non-zero, and for all $n$:
$$
\opAC(\epsilon_n,g,u_n+t_n)=0.
$$
\noindent For all $n$, denote $v_n=u_n+t_n$. Observe that the argument of Proposition \procref{PropFirstPropernessResult} is uniform in $\epsilon$ as $\epsilon$ tends to $+\infty$, and there therefore exists $v_\infty\in C^{\lambda+2}$ towards which $(v_n)_\ninn$ subconverges. For all $n$:
$$
\Delta_g v_n - \epsilon_n^{-1}f(v_n) = 0.
$$
\noindent Upon taking limits, it follows that $\Delta_g v_\infty=0$, and so $v_\infty$ is equal to a constant, $c$, say. On the other hand, for all $n$:
$$
\int f(v_n)\opdVol = \int \epsilon_n\Delta_g v_n\opdVol = 0,
$$
\noindent and upon taking limits, it follows that:
$$
f(c)\opVol(M) = \int f(c)\opdVol = 0,
$$
\noindent and so $c$ is a zero of $f$. In particular, $\Delta_g(\epsilon_n u_n)=(f(v_n))_\ninn$ converges to $0$ in the $C^{\lambda}$-topology. However, by the Closed Graph Theorem, the restriction of $\Delta_g$ to $\opKer(\Delta_g)^\perp$ is a linear isomorphism onto its image, and it follows that $(\epsilon_n \|u_n\|_{\lambda+2})_\ninn$ converges to $0$. Finally, for all $n$:
$$
\opVol(M) t_n = \int v_n\opdVol,
$$
\noindent from which it follows that $(\left|t_n-c\right|)_\ninn$ converges to $0$. It now follows from Proposition \procref{PropPreliminaryToCondensationOfSolutions} that for sufficiently large $n$, $u_n=0$. This absurd by hypothesis, and the result follows.\qed
\newsubhead{Morse Homology} We now study the Morse Homology of the Allen Cahn Equation. The construction is fairly standard, and we refer the reader to our forthcoming paper \cite{SmiMH} for a detailled outline in the H\"older space framework. We assume henceforth that $\opDim(M)\geqslant 3$. Let $g$ be as in Theorem \procref{ThmTheSolutionSpaceIsGenericallyNice} and let $\epsilon$ be such that $\epsilon^{-1}\notin\opSpec(-\Delta_g)$.
\medskip
\noindent For all $k\in\Bbb{N}$, we define $\Cal{Z}_{\epsilon,g,k}\subseteq\Cal{Z}_{\epsilon,g}$ by:
$$
\Cal{Z}_{\epsilon,g,k} = \left\{ u\in\Cal{Z}_{\epsilon,g}\ |\ \opIndex(u)=k\right\},
$$
\noindent and for all $k\in\Bbb{N}$, we define the {\bf chain group} $C_k$ by:
$$
C_k = \Bbb{Z}_2[\Cal{Z}_{\epsilon,g,k}] = \left\{ f:\Cal{Z}_{\epsilon,g,k}\rightarrow\Bbb{Z}_2\right\}.
$$
\noindent Morse Homology theory defines a canonical chain mapping $\partial_k:C_k\rightarrow C_{k-1}$ in terms of solutions to the parabolic Allen-Cahn Equation, $\oppAC_{\epsilon,g}:=\partial_t - \opAC_{\epsilon,g}$, over the space $\Bbb{R}\times M$. The Morse Homology of the Allen-Cahn Equation is then defined to be the homology of the chain complex $(C_*,\partial_*)$. That is, for all $k$:
$$
\opHAC_k = \frac{\opKer(\partial_k)}{\opIm(\partial_{k+1})}.
$$
\noindent Importantly, $\opHAC_*$ is independant, up to isomorphism, of the pair $(\epsilon,g)$ used to define it. In actual fact, the preceeding construction would require that all elements of $\Cal{Z}_{\epsilon,g}$ be non-degenerate. However, since all critical points of $e_g$ are themselves non-degenerate, we use a perturbation argument to show that degenerate elements of $\Cal{Z}_{\epsilon,g}$ do not contribute to the homology: in other words, we simply ignore them. The justification is analogous to the manner in which the function $F_\epsilon(t):=t^3 +\epsilon t$ has a degenerate critical point at $0$ when $\epsilon=0$, and no critical points for $\epsilon>0$, in contrast to the function $G_\epsilon(t)=t^4+\epsilon t^2$, which has a critical point at $0$ for all $\epsilon$.
\medskip
\noindent In order to calculate the Morse Homology, we suppose that $\epsilon\gg 0$. By Proposition \procref{PropCondensationOfSolutions}, we may suppose that $\Cal{Z}_{\epsilon,g}$ only consists of constant solutions, and furthermore, by Proposition \procref{PropConstantSolutions}, we may suppose that the Morse Index of the constant solution $u=c$ is equal to $0$ or $1$ according as $f'(c)$ is positive or negative respectively. Let $F$ be any primitive of $f$, let $c_\pm$ be zeroes of $f$, and let $w:\Bbb{R}\rightarrow\Bbb{R}$ be such that:
$$
\partial_t w = - f\circ w,\qquad \mlim_{t\rightarrow\pm\infty}=c_\pm.
$$
\noindent That is, $w$ is a gradient flow of $F$ from $c_-$ to $c_+$. We extend $w$ to a function from $\Bbb{R}\times M$ into $\Bbb{R}$ by setting it to be constant in the $x$ direction. Observe that $w$ is then a bounded solution to the parabolic Allen-Cahn Equation. That is:
$$
\oppAC_{\epsilon,g}w = (\partial_t - \opAC_{\epsilon,g})w = 0.
$$
\noindent We therefore refer to such a function $w$ as a {\bf space-constant trajectory}. As in the elliptic case, we say that $w$ is {\bf non-degenerate} whenever the linearisation of $\oppAC_{\epsilon,g}$ around $w$ defines a surjective mapping from the inhomogeneous Sobolev space $H^{1,2}(\Bbb{R}\times\Bbb{M})$ into $L^2(\Bbb{R}\times\Bbb{M})$. In order to correctly calculate the Morse Homology, we have to show that all trajectories that we study are non-degenerate. However:
\proclaim{Proposition \nextprocno}
\noindent For all $g\in\Cal{M}^{\mu+1}$, there exists $B>0$ such that for $\epsilon>B$, every space-constant trajectory is non-degenerate.
\endproclaim
\proclabel{PropConstantTrajectoriesAreNonDegenerate}
\proof Let $w:\Bbb{R}\times M\rightarrow\Bbb{R}$ be a space constant trajectory, and let $L$ be the linearisation of $\oppAC_{\epsilon,g}$ about $w$. For all $\varphi:\Bbb{R}\times M\rightarrow M$:
$$
L\varphi = (\partial_t - \epsilon\Delta_g)\varphi - (f'\circ w)(t)\varphi.
$$
\noindent By the Sturm-Liouville Theorem, there exists an orthonormal basis $(\psi_n)_\ninn$ of $L^2(M)$ consisting of eigenfunctions of $-\Delta_g$. Let $0=\lambda_0<\lambda_1\leqslant...$ be the corresponding eigenvalues. Define $B>0$ such that $B>\|f'\|_{L^\infty}/\lambda_1$. We claim that $B$ has the desired properties. Indeed, choose $\epsilon>B$. By Proposition \procref{PropConstantSolutions}, both $c_-$ and $c_+$ are non-degenerate with Morse Indices equal either to $0$ or $1$. Observe, moreover, that $\opIndex(c_-)=1$ and $\opIndex(c_+)=0$. By the Atiyah-Patodi-Singer Index Theorem (c.f. \cite{Salamon}), $L$ defines a Fredholm mapping from $H^{1,2}(\Bbb{R}\times M)$ into $L^2(\Bbb{R}\times M)$ of Fredholm index equal to $1$. Thus, in order to show that $w$ is non-degenerate, it suffices to show that $\opDim(\opKer(L))\leqslant 1$. However, choose $\varphi\in\opKer(L)$. For $k\geqslant 1$, define $\varphi_k:\Bbb{R}\rightarrow\Bbb{R}$ by $\varphi_k(t)=\langle\varphi_t,\psi_n\rangle$. Observe that $\varphi_k\in L^2(\Bbb{R})$. However:
$$
\dot{\varphi}_k = (\epsilon\lambda_k + (f'\circ w)(t))\varphi_k.
$$
\noindent Since $\epsilon>B$, there exists $\delta>0$ such that $(\epsilon\lambda_k + (f'\circ w)(t))>\delta$. Thus, over any interval in which $\varphi_n$ is non-vanishing, we have:
$$
\partial_t(\opLog(\left|\varphi_n\right|)) \geqslant \delta,
$$
\noindent and since $\psi_k\in L^2(\Bbb{R})$, it must therefore vanish identically. $\varphi_t$ therefore lies in the linear span of $\psi_0$ for all $t$. That is, it is constant in space. However, since the space of solutions to a first order ODE is at most $1$-dimensional, it follows that $\opKer(L)$ is also at most one-dimensional, and we conclude that $w$ is non-degenerate, as desired.\qed
\medskip
\noindent This allows us to calculate the Morse-Homology:
\proclaim{Proposition \nextprocno}
\noindent The Morse Homology of the Allen-Cahn Operator is given by:
$$
H\opAC_k = \left\{\matrix\Bbb{Z}_2\ \text{if}\ k=0,\hfill\cr 0\ \text{otherwise.}\hfill\cr\endmatrix\right.
$$
\endproclaim
\proof Let $B$ be as in Proposition \procref{PropConstantTrajectoriesAreNonDegenerate} and choose $\epsilon>B$. Upon increasing $B$ is necessary, it follows from Propositions \procref{PropConstantSolutions} and \procref{PropCondensationOfSolutions} that $\Cal{Z}_{\epsilon,g}$ only consists of constant solutions and, moreover, that if $u=c$ is a constant solution, then it is non-degenerate and its Morse Index is equal to $0$ or $1$ according as $f'(c)$ is positive or negative respectively. By Property \eqnref{EqnPropertiesOfPotentialFunction} of $f$, $f$ has an odd number of zeroes, $c_1<...<c_{2n+1}$. Moreover, if $k$ is odd, then $f'(c_k)>0$, and if $k$ is even, then $f'(c_k)<0$. Consequently:
$$
\Cal{Z}_{d,0} = \left\{c_1,c_3,...,c_{2n+1}\right\},\qquad \Cal{Z}_{d,1} = \left\{c_2,c_4,...,c_{2n}\right\},
$$
\noindent and $\Cal{Z}_{d,p}$ is empty for all $p\geqslant 2$. In particular, for all $p\geqslant 2$, $C_p=0$ and so $H\opAC_p=0$. For $1\leqslant k\leqslant n$, there are two space-constant trajectories leaving $c_{2k}$, terminating in $c_{2k-1}$ and $c_{2k+1}$ respectively. Moreover, by Proposition \procref{PropConstantTrajectoriesAreNonDegenerate}, these space-constant trajectories are non-degenerate, and by the unstable manifold theorem (c.f. \cite{Weber}), up to reparametrisation in time, there are no other bounded solutions $w_t(\cdot):=w(t,\cdot)$ to the parabolic Allen-Cahn Equation which converge to $c_{2k}$ as $t$ tends to minus infinity. It follows from the definition of the chain map (c.f. \cite{SmiMH}) that:
$$
\partial_1 c_{2k} = c_{2k-1}+c_{2k+1}.
$$
\noindent In particular, $\left\{\partial_1 c\ |\ c\in \Cal{Z}_{d,1}\right\}$ is a linearly independent subset of $C_0$, and so:
$$
\opDim(H\opAC_1)=\opDim(\opKer(\partial_1))=0.
$$
\noindent Finally, by the Rank-Nullity Theorem, $\opDim(\opIm(\partial_1))=n$, and so $\opDim(H\opAC_0)=1$, and the result now follows.\qed
\medskip
\noindent We now return to the specific case studied in the introduction where $f(u)=u^3-u$, and we prove Theorem \procref{ThmBifurcation}:
\medskip
{\bf\noindent Proof of Theorem \procref{ThmBifurcation}:\ } For all $k$, we define $X_k\subseteq\Cal{Z}_{\epsilon,g}$ to be the set of all stationary solutions of Morse-Index equal to $k$, and we define $C_k$ and $\partial_k$ as outlined above. Denote $l=\opIndex(0)$. Since $f$ is odd, multiplication by $-1$ maps $\Cal{Z}_{\epsilon,g}$ to itself, and all solutions of $\opAC_{\epsilon,g}u=0$ which are different to $0$ therefore exist in pairs. The set $X_k$ therefore has even cardinality for all $k\neq l$ and odd cardinality when $k=l$. In other words, $C_{k}$ has odd dimension for $k\neq l$ and even dimension for $k=l$. For all $k$, let $K_k$ be the kernel of $\partial_k$. We claim that $K_k$ is odd-dimensional for all $0<k<l$. Indeed, choose $0<k<l-1$ and suppose that $K_k$ is odd-dimensional. Then, since $\opHAC_k=0$, it follows that the image of $\partial_{k+1}$ is also odd-dimensional, and since $C_{k+1}$ is even-dimensional, it follows by the Rank-Nullity Theorem that $K_{k+1}$ is odd-dimensional. However, since $K_0=C_0$ and since $H\opAC_0=\Bbb{Z}_2$, $\opIm(\partial_1)$ is also odd-dimensional, and it follows by the Rank-Nullity Theorem that $K_1$ is also odd-dimensional. We conclude by induction that $K_k$ is odd dimensional for all $0<k<l$ as asserted. In particular, for all $0<k<l$, $K_k$ is non-trivial, and thus so too is $C_k$, from which it follows that $X_k$ is non-empty, as desired.\qed
\goodbreak
\newhead{Bibliography}
{\leftskip = 5ex \parindent = -5ex
\leavevmode\hbox to 4ex{\hfil \cite{Nirenberg}}\hskip 1ex{Agmon S., Nirenberg L., Lower bounds and uniqueness theorems for solutions of differential equations in a Hilbert space, {\sl Comm. Pure Appl. Math.}, {\bf 20}, (1967), 207--229}%
\medskip
\leavevmode\hbox to 4ex{\hfil \cite{AllenCahn}}\hskip 1ex{Allen S., Cahn J., A microscopic theory for antiphase boundary motion and its application to antiphase domain coarsening, {\it Acta Metall.}, {\bf 27}, (1979), 1084--1095}
\medskip
\leavevmode\hbox to 4ex{\hfil \cite{Aronszajn}}\hskip 1ex{Aronszajn N., A unique continuation theorem for solutions of elliptic partial differential equations or inequalities of second order, {\sl J. Math. Pures Appl.}, {\bf 36}, (1957), 235--249}%
\medskip
\leavevmode\hbox to 4ex{\hfil \cite{Cabre}}\hskip 1ex{Cabr\'e X., Uniqueness and stability of saddle-shaped solutions to the Allen-Cahn equation, {\sl J. Math. Pures Appl.}, {\bf 98}, (2012), no. 3, 239--256}
\medskip
\leavevmode\hbox to 4ex{\hfil \cite{DeGiorgi}}\hskip 1ex{De Giorgi E., Some conjectures on flow by mean curvature, {\sl White Paper}, (1990)}%
\medskip
\leavevmode\hbox to 4ex{\hfil \cite{GilbTrud}}\hskip 1ex{Gilbarg D., Trudinger N. S., {\sl Elliptic partial differential equations of second order}, Classics in Mathematics, Springer-Verlag, Berlin, (2001)}%
\medskip
\leavevmode\hbox to 4ex{\hfil \cite{Ilmanen}}\hskip 1ex{Ilmanen T., Convergence of the Allen-Cahn Equation to Brakke's motion by Mean Curvature, {\sl J. Diff. Geom.}, {\bf 38}, (1993), 417--461}
\medskip
\leavevmode\hbox to 4ex{\hfil \cite{Kato}}\hskip 1ex{Kato T., {\sl Perturbation theory for linear operators}, Classics in Mathematics, Springer-Verlag, Berlin, (1995)}%
\medskip
\leavevmode\hbox to 4ex{\hfil \cite{Palais}}\hskip 1ex{Palais R. S., Morse theory on Hilbert manifolds, {\sl Topology}, {\bf 2}, (1963), 299--340}
\medskip
\leavevmode\hbox to 4ex{\hfil \cite{Salamon}}\hskip 1ex{Robbin J., Salamon D., The spectral flow and the Maslov index, {\sl Bull. London Math. Soc.}, {\bf 27} (1995), no. 1, 1--33}%
\medskip
\leavevmode\hbox to 4ex{\hfil \cite{SalamonWeber}}\hskip 1ex{Salamon D., Weber J., Floer homology \& the heat flow, {\sl GAFA}, {\bf 16}, (2006), no. 5, 1050--1138}%
\medskip
\leavevmode\hbox to 4ex{\hfil \cite{Schwarz}}\hskip 1ex{Schwarz M., {\sl Morse homology}, Progress in Mathematics, {\bf 111}, Birkh\"auser Verlag, Basel, (1993)}%
\medskip
\leavevmode\hbox to 4ex{\hfil \cite{Smale}}\hskip 1ex{Smale S., An infinite dimensional version of Sard's theorem, {\sl Amer. J. Math.}, {\bf 87}, (1965), 861--866}%
\medskip
\leavevmode\hbox to 4ex{\hfil \cite{SmiMH}}\hskip 1ex{Smith G., A H\"older Space Approach to Morse Homology, in preparation}%
\medskip
\leavevmode\hbox to 4ex{\hfil \cite{Taylor}}\hskip 1ex{Taylor M. E., {\sl Partial differential equations III. Nonlinear equations.}, Applied Mathematical Sciences, {\bf 117}, Springer, New York, (2011)}%
\medskip
\leavevmode\hbox to 4ex{\hfil \cite{Weber}}\hskip 1ex{Weber J., Morse homology for the heat flow, {\sl Math. Z.}, {\bf 275}, (2013), no. 1, 1--54}%
\par}

%
%
%
\enddocument

%% file: preamble.tex
%
%
%
\let\myfrac=\frac%
\input eplain %
\let\frac=\myfrac%
\input amstex \input epsf %
%
%
\loadeufm\loadmsam\loadmsbm\message{symbol names}\UseAMSsymbols\message{,}%
\magnification 1200 %
\font\myfontdefault=cmr10%
\newif\ifmakebiblio%
\newif\ifinappendices%
\newif\ifundefinedreferences%
\newif\ifchangedreferences%
\makebibliofalse%
\undefinedreferencesfalse%
\changedreferencesfalse%
%
%
%
%
%
\def\setcatcodes{\catcode`\!=0 \catcode`\\=11}%
{\global\let\noe=\noexpand%
\catcode`\@=11 \catcode`\_=11 \setcatcodes%
!global!def!_@@internal@@makeref#1{%
!global!expandafter!def!csname #1ref!endcsname##1{%
!csname _@#1@##1!endcsname%
!expandafter!ifx!csname _@#1@##1!endcsname!relax%
    !write16{#1 ##1 not defined - run saving references}%
    !undefinedreferencestrue%
!fi}}%
!global!def!_@@internal@@makelabel#1{%
!global!expandafter!def!csname #1label!endcsname##1{%
!edef!temptoken{!csname #1info!endcsname}%
!ifloadreferences%
    !expandafter!ifx!csname _@#1@##1!endcsname!relax%
        !write16{#1 ##1 not hitherto defined - rerun saving references}%
        !changedreferencestrue%
    !else%
        !expandafter!ifx!csname _@#1@##1!endcsname!temptoken%
        !else%
            !write16{#1 ##1 reference has changed - rerun saving references}%
            !changedreferencestrue%
        !fi%
    !fi%
!else%
    !expandafter!edef!csname _@#1@##1!endcsname{!temptoken}%
    !edef!textoutput{!write!references{\global\def\_@#1@##1{!temptoken}}}%
    !textoutput%
!fi}}%
!global!def!makecounter#1{!_@@internal@@makelabel{#1}!_@@internal@@makeref{#1}}%
!unsetcatcodes%
}
%
%
%
%
%
\def\turnintolatin#1{\ifcase #1 _\or i\or ii\or iii\or iv\or v\or vi\or vii\or viii\or ix\or x\or xi\or xii\or xiii\or xiv\or xv\or xvi\or xvii\or xviii\or xix\or xx\or xxi\or xxii\or xxiii\or xxiv\or xxv\or xxvi\fi}%
\def\alphanum#1{\ifcase #1 _\or A\or B\or C\or D\or E\or F\or G\or H\or I\or J\or K\or L\or M\or N\or O\or P\or Q\or R\or S\or T\or U\or V\or W\or X\or Y\or Z\fi}%
\newwrite\references%
\ifloadreferences{\catcode`\@=11 \catcode`\_=11 \input references.tex }%
\else{\openout\references=references.tex }%
\fi%
%
%
\newcount\headno%
\global\headno=0%
\def\headinfo{\ifinappendices\alphanum\headno\else\the\headno\fi}%
\def\nextheadno{\global\advance\headno by 1 \global\subheadno=0 \global\procno=0 \headinfo}%
\makecounter{head}%
%
%
\newcount\subheadno%
\global\subheadno=0%
\def\subheadinfo{\headinfo.\the\subheadno}%
\def\nextsubheadno{\global\advance\subheadno by 1 \global\procno=0 \subheadinfo}%
\makecounter{subhead}%
%
%
\newcount\procno%
\global\procno=0%
\def\procinfo{\subheadinfo.\the\procno}%
\def\nextprocno{\global\advance\procno by 1 \procinfo}%
\makecounter{proc}%
%
%
\newcount\figno%
\global\figno=0%
\def\figinfo{\subheadinfo.\the\figno}%
\def\nextfigno{\global\advance\figno by 1 \figinfo}%
\makecounter{fig}%
%
%
\newcount\eqnno%
\global\eqnno=0%
\def\eqninfo{\text{(\alphanum{\the\eqnno})}}%
\def\nexteqnno{\global\advance\eqnno by 1 \eqninfo}%
\makecounter{eqn}%
%
%
%
%
%
\def\gobbleeight#1#2#3#4#5#6#7#8{}%
\newcount\citationno%
\global\citationno=0%
\def\citationinfo{\the\citationno}%
\makecounter{citation}%
\newwrite\biblio%
\def\newref#1#2{%
\def\temptext{#2}%
\edef\bibliotextoutput{\expandafter\gobbleeight\meaning\temptext}%
\global\advance\citationno by 1\citationlabel{#1}%
\ifmakebiblio%
    \edef\fileoutput{\write\biblio{\noindent\hbox to 0pt{\hss$[\the\citationno]$}\hskip 0.2em\bibliotextoutput\medskip}}%
    \fileoutput%
\fi}%
\def\cite#1{%
$[\citationref{#1}]$%
\ifmakebiblio%
    \edef\fileoutput{\write\biblio{#1}}%
    \fileoutput%
\fi%
}%
%
%
%
%
\let\mypar=\par%
\edef\Pagetitle={Blank}\headline={\hfil\Pagetitle\hfil}%
\edef\Pagefooter={Blank}\footline={\hfil\Pagefooter\hfil}%
%
%
\newcount\showpagenumflag%
\global\showpagenumflag=0 %
\def\nextoddpage%
{\newpage\ifodd\pageno%
\else\global\showpagenumflag=0 %
\null\vfil\eject%
\global\showpagenumflag=1 %
\fi}%
%
%
\font\headfont=cmb12%
\def\newhead#1%
{\ifhmode\mypar\fi%
\ifnum\headno=0 \else\goodbreak\bigskip\fi%
{\headfont\noindent\nextheadno\ - #1.}
\nobreak\medskip}%
%
%
\def\newsubhead#1%
{\ifhmode\mypar\fi%
\ifnum\subheadno=0 \else\goodbreak\medskip\fi%
{\bf\noindent\nextsubheadno\ - #1.\ }}%
%
%
\newif\ifinproclaim%
\global\inproclaimfalse%
\def\proclaim#1{%
\goodbreak\medskip
\bgroup\inproclaimtrue%
\noindent{\bf #1}%
\nobreak\medskip\sl}%
\def\noskipproclaim#1{%
\goodbreak\medskip%
\bgroup\inproclaimtrue%
\noindent{\bf #1}\nobreak\sl}%
\def\endproclaim{\mypar\egroup\nobreak\medskip\ignorespaces}%
%
%
%
\newcount\xpos\newcount\ypos
\def\makelabelgrid{%
\xpos=-5 \ypos=-5 %
\loop\ifnum\xpos<6 %
{\loop\ifnum\ypos<6 %
\def\labeltext{x}%
\ifnum\xpos=0\def\labeltext{+}\fi%
\ifnum\ypos=0\def\labeltext{+}\fi%
\placelabel[\xpos][\ypos]{\labeltext}%
\advance\ypos by 1 %
\repeat}%
\advance\xpos by 1 %
\repeat}%
\def\placelabel[#1][#2]#3{{%
\setbox10=\hbox{\raise #2cm \hbox{\hskip #1cm #3}}%
\ht10=0pt \dp10=0pt \wd10=0pt \box10}}%
%
%
%
%
\def\myitem#1{\noindent\hbox to .5cm{\hfill#1\hss}}%
%
%
%
%
%
%
%
%
%
\font\sansseriften=cmss10%
\font\sansserifseven=cmss7%
\font\sansseriffive=cmss5%
\newfam\sansseriffam%
\textfont\sansseriffam=\sansseriften%
\scriptfont\sansseriffam=\sansserifseven%
\scriptscriptfont\sansseriffam=\sansseriffive%
\def\mathsf{\fam\sansseriffam}%
%
%
%
\font\boldten=cmb10%
\font\boldseven=cmb7%
\font\boldfive=cmb5%
\newfam\mathboldfam%
\textfont\mathboldfam=\boldten%
\scriptfont\mathboldfam=\boldseven%
\scriptscriptfont\mathboldfam=\boldfive%
\def\mathbf{\fam\mathboldfam}%
%
%
%
\font\mycmmiten=cmmi10%
\font\mycmmiseven=cmmi7%
\font\mycmmifive=cmmi5%
\newfam\mycmmifam%
\textfont\mycmmifam=\mycmmiten%
\scriptfont\mycmmifam=\mycmmiseven%
\scriptscriptfont\mycmmifam=\mycmmifive%
\def\hexa#1{\ifcase #1 0\or 1\or 2\or 3\or 4\or 5\or 6\or 7\or 8\or 9\or A\or B\or C\or D\or E\or F\fi}%
\mathchardef\mathi="7\hexa\mycmmifam7B%
\mathchardef\mathj="7\hexa\mycmmifam7C%
%
%
\font\mymsbmten=msbm10 at 8pt%
\font\mymsbmseven=msbm7 at 5.6pt
\font\mymsbmfive=msbm5 at 4pt%
\newfam\mymsbmfam%
\textfont\mymsbmfam=\mymsbmten%
\scriptfont\mymsbmfam=\mymsbmseven%
\scriptscriptfont\mymsbmfam=\mymsbmfive%
\mathchardef\mybeth="7\hexa\mymsbmfam69%
\mathchardef\mygimmel="7\hexa\mymsbmfam6A%
\mathchardef\mydaleth="7\hexa\mymsbmfam6B%
%
%
%
%
\def\proof{{\noindent\bf Proof:\ }}%
\def\remark{{\noindent\bf Remark:\ }}%
\def\qed{~$\square$}%
\def\makeop#1{\global\expandafter\def\csname op#1\endcsname{{\text{#1}}}}%
\def\makeopsmall#1{\global\expandafter\def\csname op#1\endcsname{{\text{\lowercase{#1}}}}}%
%
%
\def\munion{\mathop{\cup}}%
\def\minter{\mathop{\cap}}%
%
%
\makeop{Ext}%
\makeop{Int}%
\makeop{Dist}%
\makeop{Diam}%
\makeop{Length}%
%
%
%
\def\ninn{{n\in\Bbb{N}}}%
\def\minn{{m\in\Bbb{N}}}%
\def\mlim{\mathop{{\text{Lim}}}}%
\def\mlimsup{\mathop{{\text{LimSup}}}}%
\def\mliminf{\mathop{{\text{LimInf}}}}%
\def\msup{\mathop{{\text{Sup}}}}%
%
%
%
\makeop{Dim}%
\makeop{Ker}%
\makeop{Coker}%
\makeop{Tr}%
\makeop{Adj}%
\makeop{Det}%
\makeop{End}%
\makeop{Lin}%
\makeop{Symm}%
\makeop{Mult}%
%
%
\makeop{dx}%
\makeop{dy}%
\makeop{dz}%
\makeop{dt}%
\makeop{dVol}%
\makeop{dArea}%
\makeop{Supp}%
\makeop{Hess}%
\makeop{Lip}%
%
%
\makeop{Re}%
\makeop{Im}%
\makeop{Arg}%
\makeop{Log}%
\makeop{Exp}%
%
%
\makeopsmall{Cos}%
\makeopsmall{Sin}%
\makeopsmall{Tan}%
\makeopsmall{Sec}%
\makeopsmall{Cosec}%
\makeopsmall{Cot}%
\makeopsmall{ArcCos}%
\makeopsmall{ArcSin}%
\makeopsmall{ArcTan}%
\makeopsmall{ArcSec}%
\makeopsmall{ArcCosec}%
\makeopsmall{ArcCot}%
%
%
\makeopsmall{Cosh}%
\makeopsmall{Sinh}%
\makeopsmall{Tanh}%
\makeopsmall{ArcCosh}%
\makeopsmall{ArcSinh}%
\makeopsmall{ArcTanh}%
%
%
\makeop{Vol}%
\makeop{Area}%
\makeop{Riem}%
\makeop{Ric}%
\makeop{Scal}%
\makeop{Euc}%
\makeop{Imm}%
\makeop{Emb}%
%
%
\makeop{Id}%
\makeop{Ad}%
\makeop{O}%
\makeop{SO}%
\makeop{SL}%
\makeop{GL}%
\makeop{Conf}%
\makeop{Homeo}%
\makeop{Diff}%
\makeop{Isom}%
%
%
\makeop{Ind}%
\makeop{Sig}%
\makeop{Spec}%
%
%
\makeop{Conv}%
\makeop{Max}%
\makeop{Min}%
\makeop{Mod}%
\makeop{Deg}%
\makeop{loc}%
%
%
%
%
%
%
%
%
%
%
%
%
%

%% file: references.tex
\global\def\_@citation@Nirenberg{1}
\global\def\_@citation@AllenCahn{2}
\global\def\_@citation@Aronszajn{3}
\global\def\_@citation@Cabre{4}
\global\def\_@citation@DeGiorgi{5}
\global\def\_@citation@GilbTrud{6}
\global\def\_@citation@Ilmanen{7}
\global\def\_@citation@Kato{8}
\global\def\_@citation@Palais{9}
\global\def\_@citation@Salamon{10}
\global\def\_@citation@SalamonWeber{11}
\global\def\_@citation@Schwarz{12}
\global\def\_@citation@Smale{13}
\global\def\_@citation@SmiMH{14}
\global\def\_@citation@Taylor{15}
\global\def\_@citation@Weber{16}
\global\def\_@proc@ThmTheSolutionSpaceIsGenericallyNice{1.1.1}
\global\def\_@proc@ThmBifurcation{1.1.2}
\global\def\_@eqn@EqnFirstPropertiesOfPotentialFunction{\relax \unhbox \voidb@x \hbox {(A)}}
\global\def\_@eqn@EqnPropertiesOfPotentialFunction{\relax \unhbox \voidb@x \hbox {(B)}}
\global\def\_@eqn@EqnAllenCahnEquation{\relax \unhbox \voidb@x \hbox {(C)}}
\global\def\_@proc@PropRegularityOfSolutions{2.1.1}
\global\def\_@proc@PropAPrioriBounds{2.1.2}
\global\def\_@proc@PropFirstPropernessResult{2.1.3}
\global\def\_@subhead@TheRegularSolutionSpace{2.2}
\global\def\_@eqn@EqnLinearisedAllenCahnEquation{\relax \unhbox \voidb@x \hbox {(D)}}
\global\def\_@proc@PropConstantSolutions{2.2.1}
\global\def\_@proc@PropPreservationOfParaproperness{2.2.2}
\global\def\_@proc@PropPiIsParaproper{2.2.3}
\global\def\_@proc@ThmSardSmale{2.2.4}
\global\def\_@proc@PropRegularSolutionSpaceIsSmooth{2.3.1}
\global\def\_@proc@PropPerturbedLaplacian{2.3.2}
\global\def\_@proc@PropSamplingPoints{2.3.3}
\global\def\_@proc@PropSurjectivity{2.3.4}
\global\def\_@proc@PropSmoothnessForGenericMetrics{2.3.5}
\global\def\_@proc@PropFirstPropositionConcerningKernels{2.4.1}
\global\def\_@proc@PropSecondPropositionConcerningKernels{2.4.2}
\global\def\_@proc@PropSolutionSpaceIsLocallyAGraph{2.4.3}
\global\def\_@proc@PropVariationOfNonDegenerateEigenvalue{2.4.4}
\global\def\_@proc@PropRelationBetweenLambdaAndEpsilon{2.4.5}
\global\def\_@proc@PropLambdaDefinedGlobally{2.4.6}
\global\def\_@proc@PropDerivativeOfLambdaIsNonVanishing{2.4.7}
\global\def\_@proc@PropAlmostNonDegenerateEpsilonForGenericMetrics{2.4.8}
\global\def\_@proc@PropGeometricPropertiesOfDegenerateCase{2.5.1}
\global\def\_@proc@PropSecondVersionGeometricPropertiesOfDegenerateCase{2.5.2}
\global\def\_@proc@PropRemovingGeodesicProperty{2.5.3}
\global\def\_@proc@PropDegenerateCaseDoesNotExistForGenericMetrics{2.5.4}
\global\def\_@proc@PropPreliminaryToCondensationOfSolutions{2.6.1}
\global\def\_@proc@PropCondensationOfSolutions{2.6.2}
\global\def\_@proc@PropConstantTrajectoriesAreNonDegenerate{2.7.1}